\documentclass[10pt, a4paper]{article}

\usepackage{amsmath,amssymb,url,graphicx, amsthm}
\usepackage[utf8]{inputenc}
\usepackage{mathrsfs}
\usepackage{mathtools,setspace}
\usepackage{fullpage}
\usepackage{authblk}

\usepackage{verbatim}

\usepackage{hyperref}

\usepackage{xcolor}
\hypersetup{
  colorlinks   = true, %Colours links instead of ugly boxes
  urlcolor     = {blue!90!black}, %Colour for external hyperlinks
  linkcolor    = {blue!90!black}, %Colour of internal links
  citecolor   = {red!90!black} %Colour of citations
}

\newcommand{\R}{\mathbb{R}}
\newcommand{\N}{\mathbb{N}}
\newcommand{\Z}{\mathbb{Z}}

\newcommand{\Go}{\Gamma _0}

\newcommand{\gtwo}{\Gamma _2}
\renewcommand{\d}{\mathrm{d}}

\newtheorem{thm}{Theorem}[section]
\newtheorem{prop}[thm]{Proposition}
\newtheorem{cor}[thm]{Corollary}
\newtheorem{lem}[thm]{Lemma}

\theoremstyle{definition}
\newtheorem{defi}[thm]{Definition}
\newtheorem{rmk}[thm]{Remark}

\author[1]{Viktor Bezborodov \thanks{Email: \texttt{viktor.bezborodov@uni-goettingen.de}}} 
\author[2]{
	Luca Di Persio \thanks{Email: \texttt{luca.dipersio@univr.it}}}

\affil[1]{%\emph{Bielefeld University, Faculty of Mathematics} \\
	{The University of Goettingen,  Institute for Mathematical Stochastics
}}

\affil[2]{
	{The University of Verona}}

\title{
	Spatial birth-and-death processes with a finite
	number of particles
}

\begin{document}

\maketitle

\begin{abstract}
	The aim of this work is 
		to establish essential properties 
		of spatial birth-and-death processes
		with general birth and death rates 
		on $\R^\d$.
  Spatial birth-and-death processes
  with time dependent rates
  are obtained
  as solutions to certain stochastic equations.
 The existence, uniqueness, uniqueness in law
 and the strong Markov property of unique solutions 
 are proven when the integral
 of the birth rate over $\R ^ \d$
 grows not faster than linearly with 
 the number of particles of the system.
Martingale properties of the constructed process
provide a rigorous connection 
to the heuristic generator.

We also study pathwise behavior of an aggregation model.
The probability of extinction
and the growth rate of
the number of particles 
conditioning on non-extinction are estimated.

\end{abstract}


\section{Introduction}

We consider spatial birth-and-death processes
with time dependent birth and death rates. 
At each moment of time
the system is 
represented as a finite collection 
of motionless  particles in $\R ^\d$. The particles can also be interpreted 
as individuals.
Existing particles 
may die and new particles may appear.
Each particle is characterized by its 
location.
 
 The state space of a spatial birth-and-death
 Markov process on $\R ^\d$ with finite number of particles
 is the space of finite subsets of $\R ^\d$
 \[
 \Gamma _0(\R ^\d)=\{ \eta \subset \R ^\d : |\eta| < \infty \},
\]
where $|\eta|$ is the number of points of $\eta$. $\Gamma_0: = \Gamma _0(\R ^\d)$ is also called the space of finite configurations.

Denote by $\mathscr{B}(\R ^\d)$ the Borel $\sigma$-algebra on $\R ^\d$.
 The evolution of the spatial birth-and-death process on $\R^\d$
 admits the following description. Let $\R _+  := [0, +\infty ) $. Two measurable functions characterize the 
 development in time,
 the birth rate 
 $b: \R ^\d \times \R _+ \times \Gamma _0 (\R^\d) \rightarrow [0,\infty)$
 and the death rate 
 $d: \R ^\d \times \R _+ \times \Gamma _0 (\R^\d) \rightarrow [0,\infty)$. 
 If the system is in state $\eta \in \Go$ 
 at time $t$, then the probability 
 that a new particle appears (a ``birth'') in a bounded set $B\in \mathscr{B}(\R ^\d)$
 over time interval $[t;t+ \Delta t]$ is 
 \[
   \Delta t \int\limits _{ B}b(x,t, \eta)dx + o(\Delta t),
 \]
 the probability that
a particle $x \in \eta$ is deleted from the configuration (a ``death'') over 
 time interval $[t;t+ \Delta t]$ is 
 \[
  d(x,t, \eta) \Delta t + o(\Delta t),
 \]
 and no two events happen simultaneously.  
 By an event we mean a birth or a death.
Using a slightly different terminology, we can say that
the rate at which a birth 
 occurs in $B$ is $\int_B b(x, t, \eta)dx$,
 the rate at which a particle $x\in \eta$
dies is $d(x,t,\eta)$, and no two events happen at the same time.

 Such processes, in which the birth and death rates
depend on the spatial structure of the system as opposed to
 classical $\Z _+$-valued birth-and-death processes (see e.g. 
 \cite[Page 116]{Branch2},
 \cite[Page 109]{Branch1}),
 were first
 studied by Preston \cite{Preston}. 
 A heuristic description
 similar to that above
 appeared already there. Our description 
 resembles the one in \cite{GarciaKurtz}.

We say that the rates $b$ and $d$, or
the corresponding birth-and-death process,
are time-homogeneous if $b$ and $d$ do not depend on time.
By abuse of notation we write in this case $b(x,s,\eta) = b (x , \eta)$, 
$d(x,s,\eta) = d (x , \eta)$.
The (heuristic) generator of a time-homogeneous spatial
birth-and-death process should be of the form

\begin{align} \label{intr generator}
L F (\eta) = \int\limits _{x \in \R^\d} b(x, \eta) [F(\eta \cup {x}) -F(\eta)] dx + 
\sum\limits _{x \in \eta} d(x, \eta ) (F(\eta \setminus {x}) - F(\eta)),
\end{align}
for
 $F$ in an appropriate domain, where $\eta \cup {x}$
 and $\eta \setminus {x}$ are shorthands for 
 $\eta \cup \{ x \}$ and
 $\eta \setminus \{ x \}$, respectively.

The purpose of this paper is twofold. First we would like to 
lay the groundwork for a rigorous analysis of spatial birth-and-death processes
with a finite number of particles. To this end we provide 
construction and the basic properties of the obtained process,
such as the strong Markov property, martingale properties,
and a coupling result 
	ensuring that under certain conditions
	one birth-and-death process dominates another.
 The approach of obtaining the process
as a solution to a certain stochastic equation
can be deemed an equivalent of the graphical representation for classical 
interacting particle systems, for example the contact process
or the voter model.
The similarity manifests itself in that in both cases the entire family of processes
starting at different possibly random times from different possibly random 
initial conditions and with different birth or death rates can be constructed from a single `noise' process.
Furthermore, the construction automatically provides a
coupling for the entire family.
The latter was used in \cite{shapenodeath} in the proof of a shape theorem;
see also \cite[Page 301]{Dur88}, \cite[Pages 33-34 and elsewhere]{Liggettbook2} for the role of the graphical representation
in the analysis of discrete-space  models.

Of course, the birth-and-death process with a finite number of particles
with time-homogeneous birth and death rates 
can be relatively easily constructed as a pure jump type Markov process (see e.g. \cite[Chapter 12]{KallenbergFound}).
However constructing a coupling for the entire infinite family of processes as described above
would be rather challenging in that framework.
Additionally, the stochastic equation approach also allows 
us to naturally incorporate the case of time-inhomogeneous
birth and death rates.
 Not much
attention has been given 
to spatial time-inhomogeneous birth-and-death processes in the mathematical literature yet, 
even though such temporally variant 
models have been shown to perform
better as  predictors in ecological
 models, see e.g. \cite{spatempVar, cattleGBR}.
Of particular interest are  periodic rates  reflecting seasonal changes.
In \cite{Statsim} a nearest-neighbor birth-and-death 
process is fitted to describe the movement of sand dunes. 
In \cite{TokyoRestaurants} spatial birth-and-death 
processes are used to describe
the process of openings and closures
of restaurants and stores in an area in Tokyo. 
The estimation of the birth and death rates
 is discussed  in \cite{BDM2021} 
   in various settings. We note that in \cite{BDM2021}  
the particles are allowed to move. 
In \cite{BolkPacala} (see also \cite{BolkPacala2})
the dynamics and moment equations
are investigated
for a model of biological population
which essentially is a birth and death 
process with rates 
that in our notation can be described by
\[
b(x, \eta ) = c_+ \sum\limits _{y \in \eta} a_+ (y-x),
\ \ \ 
d(x, \eta ) = \mu +  c_- \sum\limits _{y \in \eta} a_- (y-x).
\]
Here $\mu, c_+, c_- > 0$, and $a_+$ and $a_-$
are some kernels with compact support.
The interpretation is as follows:
\begin{itemize}
	\item  each individual
	reproduces independently of the others at a constant 
	rate $c_+ > 0$, and the offspring is displaced
	with a given kernel $a_+$;
	\item
	each individual dies at a constant rate $\mu > 0$;
	\item additionally, each individual dies 
	at rate governed by the kernel $a_-$
	due to competition. 
\end{itemize}

Further extensions 
related to this model
can be found in
  \cite{FournierMeleard, Eth04,  FKK12, Ecology, Fink36}.
  The works \cite{FournierMeleard, Eth04} 
   are focused on microscopic, or probabilistic, aspects,
   whereas \cite{FKK12, Ecology, Fink36}
   take more of a macroscopic and analytical approach. 
  Among exciting open problems for a continuous space birth-and-death process
are questions related to the asymptotic shape (see \cite{shapenodeath}  for a shape theorem
for  a spatial birth processes) and survival of the process started from   a 
single point configuration.

Our second aim is to give a detailed asymptotic analysis for the  aggregation model
and to demonstrate that it behaves  differently form the 
corresponding mesoscopic model \cite{Aggreg}.
We show certain fine asymptotic properties of the process, such as the 
finiteness of the total number of deaths over 
an infinite time interval and 
an exponential growth of the number of particles within a certain region.

\emph{A short literature overview}.
Garcia and Kurtz \cite{GarciaKurtz} obtained  
birth-and-death processes as solutions to 
certain 
stochastic integral equations for the case when the death rate $d \equiv 1$.
The systems treated there involves
an infinite number of particles.
In 
the earlier work \cite{Garcia} of Garcia  another approach was used:
birth-and-death processes were obtained
as projections of Poisson point processes. A
further development of the projection method
appears in \cite{GarciaKurtz2}.
Fournier and M{\'e}l{\'e}ard \cite{FournierMeleard} used
a similar equation 
for the construction of the Bolker--Pacala--Dieckmann--Law process
with finitely many particles.
Following ideas of \cite{GarciaKurtz}
and \cite{FournierMeleard},
we construct the  birth-and-death process described above
as a solution to a stochastic equation.

Holley and Stroock \cite{HolleyStroock} constructed the
spatial birth-and-death process
as a Markov family of  unique solutions to the corresponding martingale problem.
For the most part, they consider a process
contained in a bounded volume, with bounded 
birth and death rates.  They also
proved the corresponding result
for the nearest neighbor model in $\R ^1$
with an infinite number of particles.
Bezborodov et at. \cite{BadLattice}
construct and study  infinite particle
birth-and-death systems on the integer lattice with birth and death
rates satisfying some general conditions.
The approach taken in this paper somewhat resembles
that in \cite{BadLattice}, however in the continuous-space settings
the death part of the  
stochastic equation 
cannot be designed 
by assigning to each place its own independent Poisson process 
as is done in \cite{BadLattice}. 
Therefore the stochastic equation we use  differs significantly from the one in \cite{BadLattice}.

Belavkin and Kolokoltsov \cite{BelKol03} discuss,
among other things,
a general structure of a Feller semigroup
on disjoint unions of Euclidean spaces (see also references
therein for the construction of the Markov processes
with a given generator).
We note in this regard that
time-homogeneous birth-and-death processes
need not have the $C_0$-Feller property. 
Eibeck and Wagner \cite{EiWag03} 
discuss convergence of particle systems
to limiting kinetic equations.
In particular, they construct 
the stochastic process corresponding
to the particle system as a minimal jump process,
or pure jump type Markov process in 
the terminology of Kallenberg \cite{KallenbergFound}.
The jump kernel is assumed to be locally bounded.

\begin{comment}
Kondratiev and Skorokhod \cite{KondSkor}
consider a branching random walk, or a
contact process in continuum, with an infinite  number of particles.
That process  can be described 
as the spatial birth-and-death process with
\[
b(x,\eta)  = \lambda \sum\limits _{y \in \eta} a(x-y), 
\ \ \  d(x,\eta) \equiv 1,
\]
where $\lambda >0$ and $0\leq a \in L ^1 (\R ^\d)$.
Under some additional assumptions, 
they showed existence of
the process for a broad class of 
initial conditions.
\end{comment}

The scheme proposed by 
Etheridge and Kurtz \cite{EK14} covers
a wide range of interactions
and applies to discrete and continuous models.
Their approach is based on, among other things, 
assigning a certain mark (`level')
to each particle and letting this mark evolve
according to some law. 
A critical event, such as birth or death,
occurs when the level hit some threshold.
Recurrence properties of birth-and-death 
processes and convergence to the invariant distribution
are analyzed by {M{\o}ller} \cite{Mol89}.
Shcherbakov and Volkov \cite{ShcherbakovVolkol} consider 
the long term behavior of birth-and-death
processes on a finite graph with constant death rate and the 
birth rate of a special exponential form.  
Density bounds, the existence of an invariant measure, and certain return 
times are studied in \cite{fecundity}. 
A birth-and-death process  with constant birth rate 
involving infinitely many particles
was constructed in \cite{randomconnectiondeath} using a completely different approach 
based on a comparison with a Poisson random connection graph.
In \cite{BdP17} it is shown that the Lebesgue-Poisson measure
is a maximal irreducible measure.
Bezborodov et al. \cite{shapenodeath}
prove a shape theorem for a wide class of continuous-space birth processes
which match the above description with the death rate $d \equiv 0$.
The stochastic equation used in  \cite{shapenodeath}
to construct the process is a special case of our equation \eqref{se}.
Age-dependent birth-and-death processes
and their scaling limits are studied 
in \cite{age_dependent}.

%   Bezborodov \cite{Bez15}
%  obtains the spatial time-homogeneous birth-and-death process
%  as a unique solution to an equation slightly different
%  from the one we use here. Various questions
%  not treated in this paper are considered there,
%  for example the possibility of an explosion,
%  continuous dependence on initial conditions
%  and related semigroup of operators.

In the aforementioned references as well as in the present work
the system is represented by a Markov process.
An alternative approach consists in using the concept
of statistical dynamics that substitutes the notion of a Markov
stochastic process. This approach is based on considering evolutions of measures 
and their correlation functions. For details see e.g. 
\cite{FKK12}, \cite{Aggreg}, and references therein.

Finkelshtein et al.
\cite{Aggreg} consider different aspects of statistical dynamics 
for the aggregation model. In this model the death rate
is given by 
\[
d(x,\eta) = \exp \Big( - \sum\limits _{y \in \eta \setminus x} \phi (x-y) \Big),
\]
where $\phi$ is a positive measurable function. For more 
details see \cite{Aggreg}.
In this paper we present an
analysis of the
long time
behavior of a
 microscopic version
of this model. In particular, 
we estimate the probability of extinction and 
the speed of growth of the average number of particles.

The paper is organized as follows. 
Notation, definitions and results
are given in Section \ref{sec2}.
Proofs are collected in Sections \ref{sec3} and \ref{sec4}, 
with two auxiliary results located to Section \ref{sec5}.

 \section{The set-up and main results}\label{sec2}

 \subsection{Construction and basic properties}

The state space of a continuous-time, continuous-space
birth and death process with a finite number of particles is
 \[
 \Gamma _0(\R ^\d)=\{ \eta \subset \R ^\d : |\eta| < \infty \},
\]
where $|\eta|$ is the number of points of $\eta$.
$\Gamma _0(\R ^\d)$ is often called the 
\emph{space of finite configurations}.
The space of $n$-point configuration is
$\Gamma _0 ^{(n)} (\R ^\d) :=
\{ \eta \subset \R ^\d : |\eta| =n \}
\subset \Gamma _0(\R ^\d)$.
We will use $\Gamma _0$ and $\Gamma _0 ^{(n)}$ as  shorthands
for $\Gamma _0(\R ^\d)$ and $\Gamma _0 ^{(n)}(\R ^\d)$,
respectively.
  For  $\eta, \zeta \in \Go$, $|\eta|=|\zeta|>0$, 
we define
\begin{equation}\label{chirp}
  \rho(\eta,\zeta) :=   
  \min\limits _{\varsigma} \max\limits _{x\in \eta}
  \{|\varsigma(x) - x| \},
\end{equation}
where 
minimum is taken over the set of all bijections 
$\varsigma : \eta \to \zeta$. Note that
in \eqref{chirp} the notation $|\cdot|$
is used for the Euclidean distance in $\R ^\d$
(as opposed to the number of points as in $|\eta|$),
which hopefully should not lead to ambiguity.
Define a metric
 $\tilde \rho$ on $\Gamma _0(\R ^\d)$ by setting
 $\tilde \rho (\eta,\zeta) = 1 \wedge \rho(\eta,\zeta)$
 if $|\eta| = |\zeta| > 0$,
 $\tilde \rho (\varnothing,\varnothing) = 0$,
  and $\tilde \rho (\eta,\zeta) = 1$
 if $|\eta| \ne |\zeta|$. Denote by $\mathscr{B}(\Gamma _0 )$
 the Borel $\sigma$-algebra generated by $\tilde \rho$.
For  $\eta \in \Go$ and $a>0$ set 
\[
 \mathbf{B}_{\rho} (\eta , a): = \{ \zeta \in \Gamma ^{(|\eta|)}_0 \mid \rho(\eta,\zeta) \leq a \}.
\]
Note that
\[
 \mathscr{B}(\Gamma _0 ) = 
 \sigma \big( \{\varnothing \}, 
 {\mathbf{B}}_{\rho} (\eta , a), \eta \in \Go, a>0
  \big).
\]

 Let $X$ be a locally compact separable metric space 
 (in this paper $X$
 will be a subset of $\R ^m$ for some $m \in \N$).
 Even though the our solution process will stay 
 in $\Go$, we  introduce now a more general
 configuration space 
 to accommodate the
 driving process. 
 Denote by $\Gamma (X)$ the space of locally finite subsets
 of $X$
 \[
  \Gamma (X) = \{ \gamma \subset X \mid |\gamma \cap K| < \infty \text{ for all compact } K  \},
 \]
 also called \emph{the space of configurations over} $X$.
 The space $\Gamma (X)$ can be endowed with 
the $\sigma$-field $\mathscr{B}(X)$ generated by the projection 
 maps 
 \[
   \Gamma (X) \ni \gamma \mapsto |\gamma \cap B| \in \Z _+
 \]
where $B$ is an arbitrary bounded Borel subset of $X$.

\emph{Convention}.
With a slight abuse of notation,
we identify
 $\gamma \in \Gamma$ 
with the induced point measure on $X$, so that
\[
 \gamma ( B) = |\gamma \cap B|.
\]
This  convention also applies to elements of $\Gamma_0$ and other point processes
and is used throughout the paper.

For more details about the notions introduced here see e.g. \cite{DaleyVere},
\cite[Chapter 12]{KallenbergFound} or \cite{KondKuna}.
 Throughout this paper $\gtwo$ stands for
 $\Gamma((0, +\infty) \times \R _+)$.
Let $\pi$ be the distribution of a Poisson random measure 
on $(\gtwo, \mathscr{B}(\gtwo))$, with the intensity measure being 
the Lebesgue measure on $(0, +\infty) \times \R _+$
(here and throughout $\mathscr{B}(X)$
is the Borel $\sigma$-algebra of $X$).
Let $\mathscr{B} _t (\gtwo)$ be
the smallest sub-$\sigma$-algebra of $\mathscr{B}(\gtwo)$
such that for every $A _1 \in \mathscr{B} ((0,t])$,
$A _2 \in \mathscr{B} (\R _+)$ the map
\[
   \gtwo \ni \gamma \mapsto \gamma (A_1 \times A _2 ) \in \Z _+ 
   \cup \{+\infty\}
\]
is $\mathscr{B} _t(\gtwo)$-measurable.
Similarly, define $\mathscr{B}_{>t} (\gtwo)$
as 
the smallest sub-$\sigma$-algebra of $\mathscr{B}(\gtwo)$
such that for every $A _1 \in \mathscr{B} ((t, \infty))$,
$A _2 \in \mathscr{B} (\R _+)$ the map
\[
   \gtwo \ni \gamma \mapsto \gamma (A_1 \times A _2 ) \in \Z _+
   \cup \{+\infty\}
\]
is $\mathscr{B}_{>t} (\gtwo)$-measurable.

Let $\eta _0$ be a (random) finite 
initial  configuration,
and let $\hat {\eta} _0$ be the point process on 
$\R ^\d \times \gtwo $ obtained by associating 
to each point in $\eta _0$ an independent
Poisson point process on $\R _+ \times \R _+$,
with the distribution $\pi$. That is, 
if $\eta _0 = \sum\limits _{i=1} ^{|\eta _0|} \delta _{ x_i}$, then
\[
 \hat {\eta} _0 = \sum _{i=1} ^{|\eta _0|} \delta _{( x_i,\gamma _i)},
\]
where $\{\gamma _i\}$ is an independent collection of Poisson point processes on $\gtwo$.

 Let us now  introduce a stochastic
 differential equation driven by a Poisson point process
  designed 
in such a way that its solution is going to be 
 a spatial birth-and-death process with birth
 and death rates $b$ and $d$.
 It is not unusual to construct processes with jumps as
solution to certain stochastic equations \cite{GarciaKurtz, FournierMeleard, SDE_weak_sol_def}. A short discussion 
of the structure of \eqref{se} can be found in Remark \ref{rem SDE meaning}.

 Consider the stochastic equation with Poisson noise
\begin{equation} \label{se}
\begin{aligned}
\eta _t (B) = \int\limits _{(0,t] \times B \times [0, \infty ) \times \gtwo }
& I _{ [0,b(x,s,\eta _{s-} )] } (u) 
I \Big\{ \int\limits_{\substack{r \in (s,t],
\\
v \geq 0}} I_{[0,d(x,r, \eta _{r-})]}(v)\gamma (dr,dv) =0 \Big\}
N(ds,dx,du,d\gamma) \\
+ \int\limits _{B \times  \gtwo } &
I \Big\{ \int\limits _{\substack{r \in (0,t],
\\
v \geq 0}} I_{[0, d(x,r, \eta _{r-})]}(v)\gamma (dr,dv) =0 \Big\} \hat {\eta} _0 (dx, d \gamma ) ,
\end{aligned}
\end{equation}
 where
$(\eta _t)_{t \geq 0}$ is a cadlag $\Go$-valued solution
 process,
 $N$ is a  Poisson point process on 
$\R_+  \times \R^\d \times \R_+ \times \gtwo $,
the mean measure of $N$ is $ds \times dx \times du \times \pi$.
We require the
processes $N$ and $\hat {\eta} _0$ to be independent of each other. 
Equation \eqref{se} is understood in the sense that the equality
holds a.s. for every bounded
$B \in \mathscr{B} (\R ^\d) $ and $t \geq 0$.

\begin{rmk}\label{rem SDE meaning}
In the first integral on the right-hand side of \eqref{se} 
$x$ is the place and 
$s$ 
is the time of birth of a new particle. 
This particle is alive as long as 
$\int _s ^t I_{[0,d(x,r, \eta _{r-})]}(v)\gamma (dr,dv) =0$,
where $(x,s,u, \gamma) \in N$. Thus, $\gamma$ is the process
`responsible' for death.
The variable $u$ is a randomizer
controlling whether birth occurs at a given time and location.
In other words, each point of the driving
Poisson process $N$ in space-time carries an
extra mark $u \in \R _+ $
(used to decide whether 
the potential birth actually occurs)
and a further two-dimensional 
Poisson process $\gamma \in \gtwo$
(used to decide when it dies).
In the death term lies the main difference to the equation
considered by
Garcia and Kurtz \cite{GarciaKurtz}.
Adapted to our notation,
 the equation there is of the form 
\begin{equation} \label{se GK}
\begin{aligned}
\eta _t (B) = \int\limits _{(0,t] \times B \times [0, \infty ) \times [0, \infty ) }
& I _{ [0,b(x,\eta _{s-} )] } (u) I \Big\{ \int\limits _{r \in (s,t]} d(x, \eta _{r-})dv  < r \Big\}
\tilde N(ds,dx,du,dr) \\
+ \int\limits _{B \times  [0, \infty ) } &
I \Big\{ \int\limits _{r \in (0,t]} d(x, \eta _{r-})dv  < r \Big\} \tilde {\eta} _0 (dx, d r ),
\end{aligned}
\end{equation}
where $\tilde N$ is a Poisson point process on 
$\R_+  \times \R^\d \times \R_+ \times \R _+ $ with
 mean measure $ds \times dx \times du \times e ^{-r} dr$,
 and $\tilde {\eta} _0$ is obtained from $\eta _0$
 by attaching an independent unit exponential to each point.
 At first glance, \eqref{se} is more complicated than \eqref{se GK},
 since the death mechanism requires a whole Poisson random measure 
 on $[0;\infty) ^2$ instead of just one exponential 
 random variable. However,
 it is more difficult \emph{a priori} to define 
  a filtration $\{  \tilde {\mathscr {F}} _t  \} _ {t\geq 0}$
  such that 
 a solution to \eqref{se GK}, if unique, should be adapted to and  possess the Markov property
 with respect to $\{  \tilde {\mathscr {F}} _t  \} _ {t\geq 0}$.
 This makes working with martingale properties of a solution to \eqref{se GK}
 more convoluted.
\end{rmk}

\emph{Conditions on} $b$, $d$ and $\eta _0$.
 The birth rate $b$ and death rate $d$ are measurable maps
 from
 $\R ^\d \times \R _+ \times \Gamma _0 $ to $ [0,\infty)$.
 We assume that
  the birth rate $b$ satisfies 
the following conditions: sublinear growth on
the second variable in the sense 
that
\begin{equation} \label{sublinear growth for b} 
\int\limits _{\R ^\d} \sup\limits _{s > 0} {b}(x, s, \eta ) dx \leq c_1|\eta| +c_2,
\end{equation}
for some constants $c_1, c_2 >0$, and
that $b(x, \cdot, \eta)$ and $d(x, \cdot, \eta)$
are left-continuous for any $x \in \R ^\d$
and $\eta \in \Go$. 

\begin{comment}
and let $d$ satisfy (\textbf{is this condition really needed?})
\begin{equation} \label{condition on d}
\forall m \in \N : \sup\limits _{
\substack{ x\in \R^\d,|\eta| \leq m  \\ s \geq 0} }d(x, s, \eta) < \infty.
\end{equation}
\end{comment}

 We also assume that

\begin{equation} \label{condition on eta _0}
 E |\eta _0| < \infty.
\end{equation}

  \begin{rmk}
	Note that we consider a very general  death rate:
	apart from measurability,
	 $d$ is only required to be  left-continuous in the second argument.
\end{rmk}

 We say that $N$ is \emph{compatible} with a
 right-continuous complete filtration  $\{ \mathscr{F} _t \}$
  if for every $t\geq 0$
 \[
 {N([0,q] \times B \times C \times \Xi)}
 \]
  is $\mathscr{F} _t $-measurable for any 
 $ q \in [0,t], 
B \in \mathscr{B} (\R^\d), C \in \mathscr{B} (\R_+)$, and $ {\Xi \in \mathscr{B}_t (\gtwo)}$,
and  also
  \[
  N((t + q',t+q'+q''] \times B' \times C' \times \Xi')
  \] 
  is independent of $\mathscr{F} _t$
  for any 
  $ q'' > q' \geq 0, 
B' \in \mathscr{B} (\R^\d)$, ${C' \in \mathscr{B} (\R_+)}$, and ${ \Xi' \in \mathscr{B}_{>t} (\gtwo)}$.
  We say that $\hat \eta _0$ is compatible with $\{ \mathscr{F} _t \}$
  if   for every $t\geq 0$
  \[
  {\hat \eta _0 ([0,q]  \times \Xi)}
  \]
  is $\mathscr{F} _t $-measurable for any 
  $ q \in [0,t]$ and $ {\Xi \in \mathscr{B}_t (\gtwo)}$,
  and  also
  \[
  \hat \eta _0 ((t + q',t+q'+q'']  \times \Xi')
  \] 
  is independent of $\mathscr{F} _t$
  for any 
  $ q'' > q' \geq 0$ and ${ \Xi' \in \mathscr{B}_{>t} (\gtwo)}$.
  Sometimes we will use the representations
  \[
   N = \sum\limits _{q \in \mathcal{I}} \delta _{(s_q, x_q, u_q, \gamma _q)}, \ \ \ 
  \hat \eta _0 = \sum\limits _{q \in \mathcal{J}} \delta _{(x_q, \gamma _q)},
  \]
where $\mathcal{I}$ and $\mathcal{J}$ are some countable disjoint sets. 
Since $N$ and $\hat \eta _0$ are independent and the intensity measure of $N$ is non-atomic,
the following  holds a.s.: 
if $q \ne q'$, $q, q' \in \mathcal{I} \cup \mathcal{J}$, then $x _q \ne x _{q'}$.

 \begin{defi} \label{weak solution}
 A \emph{(weak) solution} of equation \eqref{se} is a triple 
 $(( \eta _t )_{t\geq 0} , N )$, $(\Omega , \mathscr{F} , P) $,
 $(\{ \mathscr {F} _t  \} _ {t\geq 0}) $, where 
\end{defi}
 
  (i) $(\Omega , \mathscr{F} , P)$ is a probability space, 
and $\{ \mathscr {F} _t  \} _ {t\geq 0}$ is an increasing, right-continuous
 and complete filtration of sub-$\sigma$-algebras of $\mathscr {F}$,

  (ii)   $N$ is a
 Poisson point process on $\R_+  \times \R^\d \times \R_+ \times \gtwo $  with  
 intensity  $ds \times dx \times du \times \pi$,

  (iii) $ \eta _0  $ is a random  $\mathscr {F} _0$-measurable
  element in $\Go$
  satisfying \eqref{condition on eta _0},

  (iv) the processes $N$ and $\hat \eta _0$ are independent,
  and are compatible with $\{ \mathscr {F} _t  \} _ {t\geq 0} $,

  (v) $( \eta _t )_{t\geq 0} $ is a cadlag $\Go$-valued process
adapted to $\{ \mathscr {F} _t  \} _ {t\geq 0} $, $\eta _t \big| _{t=0} = \eta _0$,  

  (vi) all integrals in \eqref{se} are well-defined, 
\[
E \int\limits _0 ^t ds \Big[ \int\limits _{\R ^\d} b(x,s, \eta _{s-}) dx
  + \sum\limits _{x \in \eta _{s-}} d(x,s,\eta _{s-}) \Big] < \infty, \ \ \ t > 0
  \]
   and
  
  (vii) equality \eqref{se} holds a.s. for all $t\in [0,\infty )$
  and all Borel sets $B$.

Following standard convention,
we also call just the  process $( \eta _t )_{t\geq 0} $ 
 a solution.
Note that for any solution $( \eta _t )_{t\geq 0} $  to \eqref{se}
a.s.
 \begin{equation}\label{faucet}
  \bigcup_{t \geq 0 } \eta _t \subset \{x_q\mid q\in \mathcal{I} \cup \mathcal{J} \}.
 \end{equation}

  Let 
 \begin{align}
   \mathscr{S} ^{0} _t =  \sigma \bigl\{ & \eta_0 , 
 N([0,q] \times B \times C \times \Xi) ,  \label{nonchalant} \\ &
 q \in [0,t], 
B \in \mathscr{B} (\R^\d), C \in \mathscr{B} (\R_+),  \Xi \in \mathscr{B}_t (\gtwo) \bigr\}, \notag
 \end{align}
 and
 let $\mathscr{S} _t$ be the completion of $\mathscr{S} ^{0} _t$ under $P$.
 Note that $\{ \mathscr{S} _t \}_{t\geq 0} $ 
 is a right-continuous filtration, see Section \ref{app2} in the Appendix.

  \begin{defi} \label{def strong solution}
    A solution  of \eqref{se} is called \emph{strong}
 if $( \eta _t )_{t\geq 0} $ is adapted to 
$(\mathscr{S} _t, t\geq 0)$.
  \end{defi}

\begin{rmk}
 In the definition above we considered solutions as processes indexed by 
 $ t\ \in[0,\infty)$. 
 The reformulations for 
the case  $ t \in [0,T]$, $0<T<\infty$, are straightforward. 
This remark also applies to many of the results below.
\end{rmk}
 
  For  $\sigma$-algebras $\mathscr{A}_1$ and $\mathscr{A}_2$,
 let $\mathscr{A}_1 \vee \mathscr{A}_2$ 
 be the smallest  $\sigma$-algebra 
 containing both $\mathscr{A}_1$ and $\mathscr{A}_2$.
 \begin{defi}
 We say that pathwise uniqueness holds for equation \eqref{se} and an
initial distribution $ \nu $ if, whenever the triples 
$(( \eta _t )_{t\geq 0} , N)$, $(\Omega , \mathscr{F} , P) $,
 $(\{ \mathscr {F} _t  \} _ {t\geq 0}) $ and 
 $(( \bar \eta _t )_{t\geq 0} , N)$, $(\Omega , \mathscr{F} , P) $,
 $(\{ \bar{\mathscr {F}} _t  \} _ {t\geq 0}) $ are weak solutions of \eqref{se} with 
$P \{ \eta _0 = \bar{\eta} _0 \} = 1 $ and $Law (\eta _0) = \nu $, 
and such that  $N$ is compatible 
with  $ \left\{ \mathscr {F} _t \vee \bar{\mathscr {F}} _t  \right\} _ {t\in [0,T]}  $,
we have $P \{ \eta _t = \bar{\eta} _t , t \in [0,\infty ) \} = 1$
(that is, the processes $\eta , \bar{\eta} $ are indistinguishable).

\end{defi}

\begin{defi} \label{joint uniqueness in law}
 We say that \textit{joint uniqueness in law} holds for equation \eqref{se} with an initial
distribution $\nu$ if any two (weak) solutions $((\eta_t) , N)$ and 
$((\eta_t  ^{ \prime }) , N  ^{\prime} )$ of \eqref{se},
$Law(\eta _0)= Law(  \eta _0  ^{\prime})=\nu$, have the same joint distribution:

$$Law ((\eta_t) , N)
= Law ((\eta_t  ^{\prime} ), N ^{\prime}  ).$$

\end{defi}

\begin{thm}\label{core thm}
Pathwise uniqueness, strong existence and joint uniqueness
in law hold for equation \eqref{se}.
  If $b$ and $d$ are time-homogeneous, 
  then the unique solution is a strong  Markov process,
 and the family of
 push-forward measures $\{P_{\alpha}, \alpha \in \Go \}$ defined
 in Remark \ref{encompass}
 constitutes a Markov process, or a Markov family of probability measures,
 on $D_{\Go}[0,\infty)$.
\end{thm}

We call the unique solution of \eqref{se} (or, sometimes, the corresponding 
family of measures on $D_{\Go}[0,\infty)$) a
\textit{(spatial) birth-and-death Markov process}.

 \begin{rmk}
 For time-homogeneous $b$ and $d$, 
 the transition probabilities of 
 the embedded Markov chain (see e.g. \cite[Chapter 12]{KallenbergFound})
 of the birth-and-death process
 are completely described by
 \begin{align}\label{transition probabilities}
  Q(\eta, \{ \eta \setminus \{ x \} \}) &= 
  \frac{d(x,\eta)}{(B+D)(\eta)} , \ \ \ \ \ \ \ x \in \eta, \ \ \eta \in \Go,  \\
  Q(\eta, \{ \eta \cup \{ x \}, x \in U\} )
  &= \frac{\int _{x\in U} b(x, \eta) dx}{(B+D)(\eta)}, \ \ \
  U \in \mathscr{B}(\R ^\d), \eta \in \Go, \notag
 \end{align}
 where $(B+D)(\eta) = \int _{\R ^\d} b(x,\eta) dx + \sum _{x \in  \eta } d(x, \eta) $.
  \end{rmk} 

The following two propositions establish a rigorous relation 
 between the unique solution to \eqref{se} 
 and $L$ defined by \eqref{intr generator}.
 To formulate the first of them, 
 let us consider the class $\mathscr{C} _b$ of cylindrical functions $F :
\Go \to \R _+$ with bounded increments.
We say that $F$ has bounded increments if 
\[
\sup\limits _{\eta \in \Go, x \in \R ^\d}
\big| F(\eta \cup \{x\} ) - F(\eta ) \big| < \infty.
\]
We say that $F$ is cylindrical 
if for some $R=R_F>0$
 \begin{equation*}
   F (\eta) = F(\zeta) \text{ whenever } \eta \cap \mathbf{B}(\mathbf{o}_\d , R) 
   = \zeta \cap \mathbf{B}(\mathbf{o}_\d , R),  
 \end{equation*}
where $\mathbf{B}(x , R)$ is
the closed ball of radius $R$ around $x$, 
and $\mathbf{o}_\d$ is the origin in $\R ^\d$.
 We recall that the filtration $\{ \mathscr{S} _t, t \geq 0\}$
 is introduced before Definition \ref{def strong solution}.

 \begin{prop} \label{poignant}
 Let $(\eta _t) _{t \geq 0}$ be a weak solution to \eqref{se}. Then 
 for any $F \in \mathscr{C} _b$
 the process
 \begin{equation} \label{hooky}
 \begin{split}
  F(\eta _t) - \int\limits _0 ^t 
  \Bigg\{ \int\limits _{\R ^\d} b(x,s,\eta _{s-}) 
  [ F(\eta _{s-} \cup \{x\}) -  F(\eta _{s-}) ] dx 
  \\ -
   \sum\limits _{x \in \eta _{s-}} d(x, s, \eta _{s-})
   [ F(\eta _{s-} \setminus \{x\}) -  F(\eta _{s-}) ] \Bigg\}  
  ds 
  \end{split}
 \end{equation}
is an $\{ \mathscr{S} _t, t \geq 0\}$-martingale.
In particular, the integral in \eqref{hooky} is well-defined a.s.
\end{prop}

\begin{rmk}
	Assume that all  conditions we imposed 
	on $b, d$, and $\eta _0$
	 are satisfied 
	except  \eqref{condition on eta _0}.
	Then we cannot claim that 
	 $E| \eta _t | < \infty$ for  $t \geq 0$.
	However, we would still get a unique solution on $[0,\infty)$
	satisfying all the items of Definition \ref{weak solution} except $(iii)$ and $(vi)$. 
	One way to see this is to  
	consider 
	a sequence of
	 initial conditions $\{\eta _0 ^{(m)}\}_{m \in \N}$, $\eta _0 ^{(m)} \subset \eta _0$,
	such that  a.s. $|\eta _0 ^{(m)}| \leq m$
	and 
$$ P \{ \eta _0 ^{(m)} = \eta _0  
	\text{	for sufficiently large } m
		 \} = 1.$$
	We are mostly interested in the case of a non-random 
	initial condition, therefore we do not discuss 
	the case when \eqref{condition on eta _0} is not satisfied in more detail.
\end{rmk}

\begin{rmk}
 The process 
started from at  a possibly random time $\tau $ from a possibly random configuration $ \zeta _ \tau$
can be obtained from the equation 
\begin{align} \label{se tau1}
\eta _{t + \tau} (B) = \int\limits _{(\tau,\tau + t] \times B \times [0, \infty ) \times \gtwo }
& I _{ [0,b(x,s,\eta _{s-} )] } (u) I \bigg\{ \int\limits _{\substack{r \in (s, \tau + t ],
		\\
		v \geq 0}} I_{[0,d(x,r, \eta _{r-})]}(v)\gamma (dr,dv) =0 \bigg\}
N(ds,dx,du,d\gamma) \notag \\
+ \int\limits _{B \times  \gtwo } 
I \bigg\{ \int \limits _{\substack{r \in (\tau, \tau + t ],
		\\
		v \geq 0}} & I_{[0, d(x,r, \eta _{r-})]}(v)\gamma (dr,dv) =0  \bigg\} \hat {\zeta} _{\tau} (dx, d \gamma ) 
+ \zeta _{\tau} (B),  \quad
t \geq 0.
\end{align}
This is the equation of the type \eqref{se} with
the initial condition $\zeta _{\tau}$ and
 the driving process
$\overline{N}$ being  the shift of $N$ by $\tau$ as 
defined in \eqref{N tau}.
We rely here on the strong Markov property of the driving process $N$
in the sense of Proposition \ref{Strong MP}.
Of course, $\tau$ should be an $\{\mathscr{S}_t, t\geq 0  \}$-stopping time,
and $\zeta _{\tau}  $ needs to be $\mathscr{S}_{\tau}$-measurable
as a map from $(\Omega, \mathscr{S}_{\tau})$ to $ (\Gamma_0 , \mathscr{B}(\Gamma_0))$ 
and such that $E |\zeta _{\tau}| < \infty$.
Considering different pairs $(\tau, \zeta _t)$,
we obtain a coupled family of the birth-and-death processes
as mentioned in the introduction.

\end{rmk}

 We also discuss a stochastic domination of one
 birth-and-death process by another. 
Consider two
equations of the form \eqref{se},

\begin{equation} \label{2se}
\begin{split}
\xi ^{(k)} _t (B) = \int\limits _{(0,t] \times B \times [0, \infty ) \times \gtwo }
& I _{ [0,b _k(x,s,\xi ^{(k)} _{s-} )] } (u) I \bigg\{ \int \limits _{\substack{r \in (s,  t ],
		v \geq 0}} I_{[0,d _k(x,r, \xi ^{(k)} _{r-})]}(v)\gamma (dr,dv) =0 \bigg\} 
\\
\times
N(ds,dx,du,d\gamma) 
+ \int\limits _{B \times  \gtwo } &
I \bigg\{ \int \limits _{\substack{r \in (0,  t ],
		v \geq 0}} I_{[0, d_k(x,r, \xi ^{(k)} _{r-})]}(v)\gamma (dr,dv) =0 \bigg\} \hat {\xi} ^{(k)} _0 (dx, d \gamma ) , \ \  k=1,2.
\end{split}
\end{equation}

We require the initial conditions $\xi ^{(k)}_0$ and the rates
$b_k$ to $d_k$ to satisfy the conditions imposed on $\eta _0$, $b$, and $d$. Let
$(\xi ^{(k)}_t)_{t \in [0,\infty )}$ 
be the unique strong solutions.  

\begin{prop} \label{couple}
	Assume that a.s. $\xi ^{(1)}_0 \subset \xi ^{(2)}_0$, and that
	for any two finite configurations $\eta ^1 \subset \eta ^2$,
	\begin{equation} \label{culpable}
	b_1(x,s,\eta ^1) \leq b_2 (x,s, \eta ^2), \ \ \ x\in \R^\d, s \geq 0,
	\end{equation}
	and 
	\begin{equation} \label{culpable2}
	d_1(x,s,\eta ^1) \geq d_2 (x,s, \eta ^2), \ \ \ x\in \eta ^1, s \geq 0.
	\end{equation}
	
	Then a.s.
	
	\begin{equation} \label{culprit}
	\xi ^{(1)}_t \subset \xi ^{(2)}_t, \ \ \ t \in [0,\infty ).
	\end{equation}

\end{prop}

 \subsection{Aggregation model}
 
 Here we consider a specific time-homogeneous model
 which we call an aggregation model. 
 This model has a property that
the death rate decreases as the number
of neighbors grows. We treat here
the death rate given below in
 \eqref{d aggr}, and,
 in addition to previous assumptions,
we require the birth rate to
grow linearly on the number of 
points in configuration in
the sense \eqref{b aggr}.
We prove in Proposition \ref{extinction} 
that the probability of extinction is 
 small if the
initial configuration has many points in some fixed Borel set
$\Lambda \subset \R ^\d$.
Propositions \ref{verfemen}, \ref{finite death}
and Theorem \ref{trajectory-wise asymptotics}
describe the
pathwise behavior 
of the process.

Let
\begin{equation} \label{d aggr}
 d(x,\eta) = \exp \Big\{ - \sum\limits _{y \in \eta} \varphi (x-y) \Big\},
\end{equation}
 where
 $\varphi $ is a nonnegative measurable function. 
 Our prime examples are $\varphi(z) = c>0$,
 or $\varphi(z) = c I \{z \in \widetilde{\Lambda}  \} $,
 where $\widetilde{\Lambda}$ is a Borel set such that 
 $\Lambda - \Lambda \in \widetilde \Lambda$.
Theorem \ref{core thm} ensures existence 
 and uniqueness of solutions,
 and that the unique solution 
 is a pure jump type Markov process.

 More specifically, let $\Lambda$ be a measurable
 non-empty
 subset of $\R^d$. Assume that the birth rate
 and the initial condition $\eta _0$ satisfy
 \eqref{sublinear growth for b} and \eqref{condition on eta _0},
 and, besides that, the
 inequalities
 \begin{equation}\label{b aggr}
  \int\limits _{\Lambda} b(x,\eta) dx \geq c|\eta \cap \Lambda|, \ \ \ \eta \in \Go,
 \end{equation}
 and
  \begin{equation}\label{b monotonne}
  \ \ \ \ \ \ \  
  \ \ \ \ \ \  b(x,\eta^1) \leq  b(x,\eta^2), \ \ \ \eta^1, \eta^2 \in \Go, \eta^1 \subset \eta^2,
 \end{equation}
hold
for some positive $c$. Note that $\Lambda$ is of positive Lebesgue measure by \eqref{b aggr}. We assume also that 
\begin{equation}
 \inf\limits _{x,y \in \Lambda} \varphi (x-y) \geq \log a,
\end{equation}
where $a>1$. 
 
 We say that the process $(\eta _t)_{t\geq 0}$  \textit {goes extinct}
if $\inf \{ t \geq 0: \eta _t = \varnothing \} < \infty $. This infimum is called
the \textit{time of extinction}.

 We want to show that,
 the probability of extinction decays 
 exponentially fast as the number of points of initial configuration
 inside $\Lambda$ grows.
 Also, we will give a few statements 
 describing the pace of growth of the number of points in the system.
 
\begin{prop} \label{extinction}

Let $\tilde{C} >0$. Then there exists $m _0 = m _0 (\tilde{C}) \in \N$
such that, whenever $m \geq m_0 $,
\[
 P _\alpha \big\{ (\eta _t)_{t\geq 0} \text{ goes extinct }  \big\} \leq \tilde{C}^{-m}
\]
for all $\alpha$ satisfying $|\alpha \cap \Lambda| = m$.

\end{prop}

\begin{prop}\label{verfemen}
 For all $\alpha \in \Go$,
\begin{equation} \label{innately}
   P_{\alpha} \bigg( \{|\eta _t \cap \Lambda| \to \infty \} 
   \cup \{ \exists t': \forall t\geq t', |\eta _t \cap \Lambda| = \varnothing \} \bigg) =1.
\end{equation}
\end{prop}
 Note that we do not require $b(\cdot, \varnothing) \equiv 0$;
if $\int _{\Lambda} b(x,\varnothing) dx >0$, then 
\eqref{innately} implies
 \[
   P_{\alpha}  \{|\eta _t \cap \Lambda| \to \infty \}  =1.
 \]
 The next proposition is a consequence of the exponentially fast decay
of the death rate.

\begin{prop} \label{finite death}
 With probability $1$ only a finite number of deaths inside $\Lambda $ occur:
\[
P _{\alpha}\bigg\{ |\eta _t \cap \Lambda| - |\eta _{t-} \cap \Lambda| = -1 
\textrm{\ for infinitely many different \ } t \geq 0   \bigg\} =0, \ \alpha \in \Go.
\]
\end{prop}

 \begin{thm} \label{trajectory-wise asymptotics}
 
Let $\alpha \in \Go$.
 For $P_{\alpha}$-almost all $\omega \in F:= 
 \{ \lim\limits _{t \to \infty} |\eta _t \cap \Lambda|=\infty \}$ we have
\begin{equation} \label{implore}
 \liminf\limits _{t\to \infty} \frac{|\eta _t \cap \Lambda|}{e^{ct}} > 0 .
\end{equation}

\end{thm}
 
 \begin{cor} \label{growth in average}
 For all configurations $ \alpha $ with $ \alpha \cap \Lambda \ne \varnothing $,
 
 \begin{equation} \label{growth in aver}
\liminf\limits _{t \to \infty}  \frac{E _{\alpha} |\eta _t \cap \Lambda|}{e^{ct}} > 0.
  \end{equation}
\end{cor}

 \textbf{Remark}. If $\Lambda$ has a finite volume and
  the birth rate is given 
 constant within $\Lambda$, that is 
 \[
 b(x, \eta) = c_3 >0, \quad x \in \Lambda,
 \]
 then from the proofs we can conclude that Theorem \ref{trajectory-wise asymptotics}
 still holds provided that 
 we replace \eqref{implore} by 
 \begin{equation} 
 \liminf\limits _{t\to \infty} \frac{|\eta _t \cap \Lambda|}{  t} > 0 .
\end{equation}

 These two growth estimates stand in contrast to the mesoscopic behavior
 of the system \cite{Aggreg}. Theorem 5.3 in \cite{Aggreg} says
 that for some values of parameters the solution to the 
 mesoscopic equation started from sufficiently small initial 
 condition stays bounded. On the contrary, 
 the microscopic system grows whenever it survives,
 and the density always grows.

 \section{Proof of Theorem \ref{core thm} and Proposition \ref{poignant}} \label{sec3}

 Let us start with the equation
\begin{equation} \label{se pure birth}
\overline{\eta} _t (B) = \int\limits _{(0,t] \times B \times [0, \infty ) \times \gtwo }
I _{ [0, \overline {b}(x,s,\eta _{s-} )] } (u) N(ds,dx,du,d\gamma) + \eta _0 (B),
\end{equation}
where
$\overline{b}(x, \eta ) := \sup\limits _{s >0, \xi \subset \eta} b(x, s, \xi) $.
Note that $\overline{b}$ satisfies sublinear growth condition
\eqref{sublinear growth for b} if $b$ does.

 This equation
is of the type \eqref{se}, 
with $\overline b$ being the birth rate
and the zero function being the death rate, and all definitions 
of existence and uniqueness of solutions
are applicable here. Later a unique solution
of \eqref{se pure birth} will be
used as
 a dominating process  to a solution to \eqref{se}.

 \begin{prop} \label{pure birth prop}
Under assumptions \eqref{sublinear growth for b} and \eqref{condition on eta _0},
strong existence and pathwise uniqueness hold for equation \eqref{se pure birth}.
In particular, the unique solution $(\bar \eta _t)_{t\geq 0}$ satisfies 
\begin{equation}\label{beschwoeren}
 E|\bar \eta _t | < \infty, \ \ \  t \geq 0.
\end{equation}
\end{prop}

 \textbf{Proof}. 
For $\omega \in  \{ \int\limits _{\R^\d }
\overline{b}(x,  \eta _0) dx =0 \} $, set $\zeta _t \equiv \eta _0$,
$\sigma _n = \infty$, $n \in \N$.

 For $\omega \in F:= \{ \int\limits _{\R^\d }
\overline{b}(x,  \eta _0) dx >0 \}$,
 we define the sequence of random pairs $\{(\sigma _n, \zeta _{\sigma _n}) \}$, where
\[
\sigma _{n+1}= \inf \bigg\{ t>0 : \int\limits _{(\sigma _n,\sigma _n+t] \times B \times [0, \infty ) \times \gtwo }
I_{[0, \overline b(x,  \zeta _{\sigma _n})]} (u) N(ds,dx,du,d\gamma) >0 \bigg\}+ \sigma _n, \ \ \sigma _0 = 0,
\]
and
\[
 \zeta _{0} = \eta _0, \ \ \ 
\zeta _{\sigma _{n+1}} = \zeta _{\sigma _n} \cup \{ z_{n+1} \}
 \]
for $z_{n+1} = \{x\in \R^\d: N (\{\sigma _{n+1}\} \times \{x \} \times [0, \overline b(x,  \zeta _{\sigma _n})] \times \gtwo
) >0  \}$.
The positions 
$z_n$ are uniquely determined almost surely on $F$. 
Furthermore, 
 $\sigma _{n+1} > \sigma _n$ a.s., and $\sigma _n$
 are finite a.s. on $F$
 (in particular because $\overline{b}(x,\zeta _{\sigma _n}) \geq \overline{b}(x,\eta _0)$).
For $\omega \in F$, we define $\zeta _t = \zeta _{\sigma _n}$
for $t\in [\sigma _n, \sigma _{n+1})$. Then by induction on $n$ 
it follows that $\sigma _n $ is a stopping time for each $n \in \N$, and 
$\zeta _{\sigma _n}$ is $\mathscr{F} _{\sigma _n} \cap F $-measurable. 
By direct substitution we see that
 $(\zeta _t)_{t\geq 0}$ is a strong
 solution to \eqref{se pure birth} on the time interval
$t\in [0, \lim\limits _{n\to \infty} \sigma _n)$. 
We are going to show that 
\begin{equation}\label{staunch}
\lim\limits _{n\to \infty} \sigma _n = \infty \ \ \ \textrm{a.s.}
\end{equation}
This relation is evidently true on the complement of $F$.
If $P(F)=0$, then \eqref{staunch} is proven.

If $P(F)>0$, define a probability measure on $F$, 
$Q(A) = \frac{P(A)}{P(F)}$, $A \in \mathscr{I} := \mathscr{F} \cap F$,
and define $ \mathscr{I}_t =\mathscr{F}_t \cap F $.

The process $N$ is independent of $F$, therefore it is a
Poisson point process on the probability space $(F,\mathscr{I}, Q)$
with the same intensity,
compatible with $\{ \mathscr{I}_t \}_{t\geq 0}$.
From now on and until other is specified, we work on the filtered
probability space $(F,\mathscr{I},\{ \mathscr{I}_t \}_{t\geq 0}, Q)$.
We use the same symbols for random processes and random variables, 
having in mind that 
we consider their restrictions to $F$.

The process $(\zeta _t)_{t\in [0, \lim\limits _{n\to \infty} \sigma _n)}$
has the Markov property, because 
the process $N$ has the strong Markov property and independent increments by Proposition \eqref{Strong MP} in the Appendix. 
Recall that for $\eta \in \Go$
and $x \in \R ^\d$, $\eta \cup x$ is a shorthand for $\eta \cup \{x\}$.
Indeed, conditioning on $ \mathscr{I}_{\sigma_n}$,
\[
 E \bigl[  I_{\{\zeta _{\sigma _{n+1}} = \zeta _{\sigma _{n}} \cup x 
 \text{ for some } x \in B \}} \mid \mathscr{I}_{\sigma_n}\bigr]=
 \frac{\int\limits _{ B} \overline b(x, \zeta _{\sigma _{n}}) dx}{\int\limits _{\R^\d }
\overline{b}(x,\zeta _{\sigma _n}) dx },
\]
thus the chain $\{ \zeta _{\sigma _{n}}\}_{n \in \Z_+}$ is a Markov chain,
and, given $\{ \zeta _{\sigma _{n}}\}_{n \in \Z_+}$,
$\sigma _{n+1} - \sigma _n $ are distributed exponentially:
\[
E\{ I_{\{\sigma _{n+1} - \sigma _n >a\}}  \mid \{ \zeta _{\sigma _{n}}\}_{n \in \Z_+}\} 
= \exp \Big\{ - a \int\limits _{\R^\d }
\overline{b}(x,\zeta _{\sigma _n}) dx \Big\}.
\]
Therefore, the random variables $\gamma _n = (\sigma _{n} - \sigma _{n-1}){(\int\limits _{\R^\d }
\overline{b}(x,\zeta _{\sigma _n}) dx)}$ 
constitute a 
sequence of independent 
 random variables exponentially distributed
with parameter $1$, independent 
of $\{ \zeta _{\sigma _{n}}\}_{n \in \Z_+}$. 
Thus
 Theorem 12.18 in \cite{KallenbergFound} implies  that
$(\zeta _t)_{t\in [0, \lim\limits _{n\to \infty} \sigma _n)}$
is a pure jump type Markov process.

The jump rate of $(\zeta _t)_{t\in [0, \lim\limits _{n\to \infty} \sigma _n)}$
is given by
\[
 c(\alpha) = \int\limits _{\R^\d }
\overline{b}(x,\alpha) dx .
\]
Condition \eqref{sublinear growth for b} implies that 
$c(\alpha) \leq c_1 |\alpha| + c_2$. Hence
\[
c (\zeta _{\sigma _n} ) \leq c_1 |\zeta _{\sigma _n}| + c_2 = 
c_1 |\zeta _0| +c_1 n  + c_2.
\]

We see that $\sum _n \frac{1}{c (\zeta _{\sigma _n} )} = \infty$ a.s., 
hence Proposition 12.19 in \cite{KallenbergFound} 
implies that $\sigma _n \to \infty$.

 Now we return again to our initial probability space 
$(\Omega,\mathscr{F},\{ \mathscr{F}_t \}_{t\geq 0}, P)$.
We have proved the existence of a strong solution. 
The uniqueness follows by induction on jumps of the process. 
Namely, let $( \tilde{\zeta} _t )_{t\geq 0} $ be another solution of \eqref{se pure birth}.
Since a.s.
$$
\int\limits _{ (0,\sigma _1) \times \R ^\d  \times [0, \infty ) \times \gtwo }
I _{ [0,0] }(u) N(ds,dx,du,d \gamma) = 0 ,
$$
(here $I _{ [0,0] }(u) = I \{ u = 0\}$)
we have $\zeta _t = \tilde \zeta _t  = \eta _ 0$ a.s. on the complement $F ^c$ for all $t\geq 0$.
 From (vii)
 of Definition \ref{weak solution} and the equality

$$
\int\limits _{ (0,\sigma _1) \times \R ^\d  \times [0, \infty ) \times \gtwo }
I _{ [0,\overline{b}(x, {\eta} _0 )] }(u) N(ds,dx,du,d \gamma) = 0 ,
$$
it follows
that $P \big( \{ \tilde{\zeta} \text{\ has a birth before \ }
\sigma _1  \} \cap F \big) = 0$. At the same time, the equality

$$
\int\limits _{ \{\sigma _1\} \times \R ^\d  \times [0, \infty ) \times \gtwo }
I _{ [0,\overline{b}(x, {\eta} _0 )] }(u) N(ds,dx,du,d \gamma) = 1 ,
$$
which holds a.s. on $F$, 
yields that $\tilde{\zeta}$ has a birth at the moment $\sigma _1$, and in the same point 
of space at that. Therefore, $\tilde{\zeta}$ coincides with $\zeta$ on $[0, \sigma _1]$ a.s. on $F$. Similar reasoning shows
that they coincide up to $\sigma _n$ a.s. on $F$, and, since $\sigma _n \to \infty$ a.s. on $F$, 

$$P \{ \tilde{\zeta} _t = {\zeta} _t \text{\ for all \ } t\geq 0 \} = 1.$$

 Thus, pathwise uniqueness holds.
 
Now we turn our attention to \eqref{beschwoeren}. Since $\zeta _t \equiv \eta _0$ on $F ^c$,
we can assume without loss of generality that $P(F) = 1$.
We can write 

\begin{gather}\label{gibberish}
|{\zeta} _t| = 
 |{\eta} _0|+  \sum\limits _{n=1}^{\infty} I \{ |\zeta _t| - |\eta _0| \geq n \} \notag
 \\  = 
  |{\eta} _0|+  \sum\limits _{n=1}^{\infty} I \{ \sigma _n \leq t \}.
\end{gather}

Since $\sigma _n = \sum\limits _{i=1}^n \frac{\gamma _i}{\int\limits _{\R^\d }
\overline{b}(x,\zeta _{\sigma _i}) dx}$, we have

\[
 \{ \sigma _n \leq t \} = \{ \sum\limits _{i=1}^n \frac{\gamma _i}{\int\limits _{\R^\d }
\overline{b}(x,\zeta _{\sigma _i}) dx} \leq t \} \subset
\{ \sum\limits _{i=1}^n \frac{\gamma _i}{c_1|\zeta _{\sigma _i}| + c_2} \leq t \}
\]
\[
\subset \{ \sum\limits _{i=1}^n \frac{\gamma _i}{(c_1+c_2)(|\eta _0|+i)} \leq t \} =
\{Z _t -Z_0 \geq n \},
\]
where
$(Z_t)$ is the Yule process,
i.e. the birth process on $\Z _+$
with transition rates 
$p_{k, k+1} =(c_1 + c_2)k$,
$p_{k, l} = 0$, $l \ne k +1$,
see, e.g., \cite[Chapter 3, Section 5]{Branch1}.
Here $(Z_t)$ is defined as follows:
$Z_t - Z _0 = n$ when 
\[
\sum\limits _{i=1}^n \frac{\gamma _i}{(c_1+c_2)(|\eta _0|+i)} \leq t 
< \sum\limits _{i=1}^{n+1} \frac{\gamma _i}{(c_1+c_2)(|\eta _0|+i)},
\]
and 
$Z_0 = |\eta _0|$. Thus, we have $| \zeta _t| \leq Z _t$ a.s.,
hence $E| \zeta _t| \leq  E Z _t < \infty$.
 The constructed
 solution is strong.
\qed

\begin{prop}\label{ex un} 
Under assumptions \eqref{sublinear growth for b}-\eqref{condition on eta _0},
 pathwise uniqueness and strong existence hold for equation \eqref{se}. 
 The unique solution $(\eta _t)$ satisfies 
 \begin{equation}\label{unfathomable}
 E| \eta _t | < \infty, \ \ \  t \geq 0.
\end{equation}
\end{prop}

\textbf{Proof}. Let us define stopping times with respect to
 ${\{ \mathscr{F} _t , t\geq 0 \}} $,
$0= \theta _0 \leq \theta _1 \leq \theta _2 \leq \theta _3 \leq ...$,
and the sequence of (random) configurations 
$\{ \eta _{\theta _j} \} _{j \in \N}$
as follows: as long
as 
\[ 
\theta _{n+1} = \theta ^b _{n+1} \wedge \theta ^\d _{n+1} + \theta _n < \infty,
\]
where 
\[
\theta ^b _{n+1}= \inf \bigg\{ t>0 : \int\limits _{ (\theta _n, \theta _n +t] \times \R^\d \times [0,\infty) \times \gtwo}
I_{[0, b(x,s,\eta _{\theta _n})]} (u) N(ds,dx,du, d \gamma) >0 \bigg\},
\]
\[
\theta ^\d _{n+1}= \inf \bigg\{ t>0 :\sum\limits _{
\substack{ q \in \mathcal{I} \cup \mathcal {J},  \\ x_q \in \eta _{\theta _n}} 
}
\int  _{(\theta _n, \theta _n +t]
 \times [0, \infty ) } 
I _{ [0,d(x_q,r,\eta _{\theta _n} )] } (v) \gamma _q (dr,dv) > 0 \bigg\},
\]
we set
${\eta} _{\theta _{n+1}}
=  {\eta} _{\theta _n} \cup \{ z_{n+1} \}$ if $\theta ^b _{n+1} \leq \theta ^\d _{n+1}$,
where
$\{z_{n+1}\} = \{z\in \R^\d: N ( \{\theta _n + \theta ^b _{n+1}\} \times \{ z \} \times \R_+ \times \gtwo) >0  \}$; 
${\eta} _{\theta _{n+1}} =  {\eta} _{\theta _n} \setminus \{ z_{n+1} \}$
if $\theta ^b _{n+1} > \theta ^\d _{n+1}$, where
$\{z_{n+1}\} = \{x_q \in {\eta} _{\theta _n}: \gamma _q ( \{\theta _n + \theta ^\d _{n+1}\} \times \R_+) >0  \}$;
the configuration $\eta _{\theta _0} = \eta _0$
is the initial condition of \eqref{se}, ${\eta} _t = {\eta} _{\theta _n} $ for 
$t \in [\theta _n , \theta _{n+1} )$.
 Note that 
\[
P\{  \theta ^b _{n+1} = \theta ^\d _{n+1}  \mid
 \min \{ \theta ^b _{n+1} , \theta ^\d _{n+1} \} < \infty \} = 0,
\]
 the points $z_n$ are a.s. uniquely determined, and 
 \[
 P\{  z_{n+1} \in {\eta} _{\theta _n}\mid \theta ^b _{n+1} \leq \theta ^\d _{n+1} \} =0.
 \]
If for some $n$
\[
 \theta _{n+1} = \infty,
\]
 we set $\theta _{n+k} =\infty$, $k\in \N$, and
$\eta _{t} =  {\eta} _{\theta _{n}}$, $t \geq {\theta _{n}}$.

Random variables $\theta _n, n \in \N$, are stopping times with respect to
the filtration ${\{ \mathscr{F} _t , t\geq 0 \}} $. 
By the strong
Markov property of a Poisson point process, see Proposition \ref{Strong MP}, we obtain that
a.s. on $\{ \theta _n < \infty \} $ the conditional distribution of $\theta ^b _{n+1}$ given
$\mathscr {F} _{\theta _n}$ is
\[
 P\left\{ \theta ^b _{n+1} > p \mid  \mathscr {F} _{\theta _n}  \right\} = 
 \exp \bigg\{ - \int_{ \theta _n} ^{\theta _n + p}ds
 \int\limits _{\R ^\d}
  b(x,s,\eta _{\theta _n}) dx  \bigg\},
\]
and a.s. on $\{ \theta _n < \infty \} $ the conditional distribution of $\theta ^\d _{n+1}$ given
$\mathscr {F} _{\theta _n}$ is 
\[
 P\left\{ \theta ^\d _{n+1} > p \mid  \mathscr {F} _{\theta _n}  \right\}
  = \exp \bigg\{ - \int_{\theta _n} ^{\theta _n + p}  ds
  \sum\limits _{x \in \eta _{\theta _n}} d(x,s,\eta _{\theta _n})
   \bigg\}.
\]
In particular,
$\theta ^b _n, \theta ^\d_n > 0$, $n \in \N$.

 We are going to show that a.s.
 \begin{equation}\label{kowtow}
  \theta _n \to \infty, \quad n \to \infty.
 \end{equation}
Denote by $\theta _k '$ the moment of the $k$-th birth. It is sufficient to show that
$\theta _k ' \to \infty$, $k\to \infty$, 
because only finitely many deaths may occur between any two births, since there
are only finitely particles. By induction
on $k '$ we can see that 
$\{ \theta _k' \} _{k' \in \N} \subset \{ \sigma _i \} _{i \in \N}$, where
$\sigma _i$ are the moments of births of $(\overline{\eta} _t)_{t\geq 0}$,
the solution of \eqref{se pure birth},
 and $\eta _t \subset \overline{\eta} _t$
 for all $t \in [0,\lim_n \theta _n)$.
For instance, let us show that $(\overline{\eta} _t)_{t\geq 0}$ has a birth at
 $\theta _1 '  $. We have  $\overline{\eta} _{\theta _1 ' -} \supset \overline{\eta} _{0} = \eta _0$, and
$\eta _{\theta _1 ' -} \subset \eta _{t} \mid _{t=0} = \eta _0 $, 
hence for all $x \in \R^\d$

$$\overline{b} (x, \overline{\eta} _{\theta _1 ' -} ) 
\geq \overline{b} (x, \eta _{\theta _1 ' -} ) \geq b (x, \theta _1 '  , \eta _{\theta _1 ' -} )$$

The latter implies that at time moment $\theta _1 '  $ a birth
 occurs for the process $(\overline{\eta} _t)_{t\geq 0}$ in the same point.
 Hence, $\eta _{\theta _1 '} \subset \overline{\eta} _{\theta _1 ' } $,
 and we can go on.
Since $\sigma _k \to \infty$ as
$k \to \infty$, we also have $\theta _k ' \to \infty$, and therefore $ \theta _n \to \infty$, $n\to \infty$. 

Let us now prove the inequality from item (vi) of Definition \ref{weak solution},
\begin{equation}\label{surreptitious}
E \int\limits _0 ^t ds \Big[ \int\limits _{\R ^\d} b(x,s, \eta _{s-}) dx
  + \sum\limits _{x \in \eta _{s-}} d(x,s,\eta _{s-}) \Big] < \infty, \ \ \ t > 0.
\end{equation}
  Denote the number of births and deaths 
before $t$ by $b_t $ and $d_t$ respectively, i.e.
\[
 b_t  = \# \{ s:| \eta _{s}| - |\eta _{s-}| = 1 \}=
\int\limits _{(0,t] \times \R ^\d \times [0, \infty ) \times \gtwo }
 I _{ [0,b(x,s,\eta _{s-} )] } (u)
N(ds,dx,du,d\gamma)
\]
and 
\[
 d_t  = \# \{ s:| \eta _{s}| - |\eta _{s-}| = - 1 \}=
\int  _{(0, t]
 \times [0, \infty ) } 
 \sum\limits _{
\substack{ q \in \mathcal{I} \cup \mathcal {J},  \\ x_q \in \eta _{r-}} 
}
I _{ [0,d(x_q,r,\eta _{r-} )] } (v) \gamma _q (dr,dv).
\]

Note that $|\eta _t| = b_t - d_t + |\eta _0|$ and $\theta _k$
are the moments of jumps for $c_t : = b_t + d_t$, 
so that
 \[
   c_t= \sum\limits _{k \in \N} I \{ \theta _k \leq t \}, \ \ \ \ t \geq 0.
  \]
For $n \in \N$ define 
\[
c ^{(n)} _t = 
\int\limits _{(0,t] \times \R ^\d \times [0, \infty ) \times \gtwo }
 I _{ [0,b(x,s,\eta _{s-} ) \wedge n] }  (u) I \{ |x| \leq n \}
N(ds,dx,du,d\gamma)
\]
\[
 +
\int  _{(0, t]
 \times [0, \infty ) }
 \sum\limits _{
\substack{ q \in \mathcal{I} \cup \mathcal {J}  \\ x_q \in \eta _{r-}} 
}
I _{ [0,d(x_q,r,\eta _{r-} ) \wedge n ] } (v) I \{ |x| \leq n \} \gamma _q (dr,dv).
\]
  Then 
  \[
   M^{(n)}_t  = c ^{(n)} _t  - \int\limits _0 ^t \int _{x:  |x| \leq n }
\big(b(x,s,\eta _{s-} ) \wedge n \big) dx ds - \int\limits _0 ^t  \sum\limits _{x \in \eta _{s-}, |x| \leq n  }
\big(d(x,s,\eta _{s-} ) \wedge n \big) ds
  \]
is a martingale with respect to
$\{ \mathscr{S} _t \}$, see e.g. \cite[(3.8), Section 3, Chapter 2]{IkedaWat}.
By the optional stopping theorem $EM^{(n)}_{\theta _1 \wedge t} = 0$,
hence 
\[
E
\int\limits _0 ^{\theta _1 \wedge t } \Bigg( \int _{x:  |x| \leq n }
b(x,s,\eta _{s-} ) \wedge n \ dx +
 \sum\limits _{x \in \eta _{s-}, |x| \leq n  }
d(x,s,\eta _{s-} ) \wedge n  \Bigg) ds = E c ^{(n)} _{t \wedge \theta _1}
 \leq  P\{\theta _{1} < t  \} \leq 1.
\]
Similarly, 
\[
E
\int\limits _{\theta _m \wedge t} ^{\theta _{m+1} \wedge t} \Bigg( \int _{x:  |x| \leq n }
b(x,s,\eta _{s-} ) \wedge n \ dx +
\sum\limits _{x \in \eta _{s-}, |x| \leq n  }
d(x,s,\eta _{s-} ) \wedge n \Bigg)  ds 
\]
\[
=E c ^{(n)} _{t \wedge \theta _{m+1}} - E c ^{(n)} _{t \wedge \theta _m}
\leq  P\{\theta _{m + 1} < t  \}.
\]
Consequently,
\[
E
\int\limits _{0} ^{t} \Bigg( \int _{x:  |x| \leq n }
b(x,s,\eta _{s-} ) \wedge n \ dx +
\sum\limits _{x \in \eta _{s-}, |x| \leq n  }
d(x,s,\eta _{s-} ) \wedge n \Bigg)  ds 
\]
\[
 \leq \sum\limits _{m=1} ^{\infty} E
 \int\limits _{\theta _m \wedge t} ^{\theta _{m+1} \wedge t} \Bigg( \int _{x:  |x| \leq n }
 b(x,s,\eta _{s-} ) \wedge n \ dx +
 \sum\limits _{x \in \eta _{s-}, |x| \leq n  }
 d(x,s,\eta _{s-} ) \wedge n \Bigg)  ds 
\]
\[
 \leq \sum\limits _{m = 1} ^ \infty  P\{\theta _{m} \leq t  \} = 
 \sum\limits _{m = 1} ^ \infty  P\{c_t  \geq m  \} = E c_t .
\]
Letting 
$n \to \infty$, we get by the monotone convergence theorem
\[
E
\int\limits _{0} ^{t} \Bigg( \int _{x \in \R ^\d }
b(x,s,\eta _{s-} )  dx +
\sum\limits _{x \in \eta _{s-} }
d(x,s,\eta _{s-} ) \Bigg)  ds  \leq E c_t .
\]
Only existing particles may disappear, hence the number of deaths $d_t $
satisfies 
\[
 d_t  \leq b_t  + |\eta _0 |.
\]
Thus, 
\begin{equation} \label{on the skids}
 E c_t  \leq 2 E b_t + E |\eta _0| 
\leq 
2 E |\bar \eta _t| + E |\eta _0| < \infty,
\end{equation}
and \eqref{surreptitious} follows.

 Since $\eta _t \subset \overline{\eta} _t$ a.s., Proposition \ref{pure birth prop} 
 implies \eqref{unfathomable}.

If follows from the above construction, \eqref{kowtow}, and \eqref{surreptitious}
that $(\eta _t)$ is a strong solution to \eqref{se}.
Similarly to the proof of Proposition \ref{pure birth prop}, 
we can show by induction on $n$ that equation \eqref{se} 
has a unique solution on $[0, \theta _n]$. Namely, each two
solutions coincide on $[0, \theta _n]$ a.s.
Thus, any solution coincides with $({\eta} _t)$ a.s. for all 
$t  \in [0, \theta _n] $.

 \qed

 \begin{rmk}\label{encompass}
 
 Assume that $b$ and $d$ are time-homogeneous.
Let $\eta _0$ be a non-random initial condition, $\eta _0 \equiv \alpha$,
$\alpha \in \Go$.
The solution of \eqref{se} with $\eta _0 \equiv \alpha$
will be denoted as $(\eta(\alpha, t))_{t\geq 0}$.
Let $P_\alpha$ be
the push-forward of $P$ under the mapping

\begin{equation}
  \Omega \ni \omega \mapsto (\eta(\alpha, \cdot)) \in D_{\Go}[0,\infty ).
\end{equation}

 It can be derived from the proof of Proposition \ref{ex un}
 that, for fixed $\omega \in \Omega$, 
the unique solution is jointly measurable in $(t, \alpha)$.
Thus, the family $\{ P_\alpha \}$
 of probability measures on $D_{\Go}[0,\infty )$
is measurable in $\alpha$,
that is, for any Borel set $\mathcal{D} \subset D_{\Go}[0,\infty ) $
the map $\Gamma_0 \ni \alpha \mapsto P_\alpha (\mathcal{D})$ is measurable.
 We will often use formulations
related to the probability space
$(D_{\Go}[0,\infty ),\mathscr{B}(D_{\Go}[0,\infty )), P_{\alpha})$;
in this case, coordinate mappings will be denoted 
by $\eta_t$,

\[
 \eta _t (x) = x(t), \ \ \   x \in D_{\Go}[0,\infty ).
\]

The processes $(\eta _t) _{t \in  [0,\infty )}$
and $(\eta(\alpha, \cdot))_{t \in  [0,\infty )} $
have the same law (under $P_{\alpha} $ and $P$, respectively).
As one would expect, the family
of measures $\{P_{\alpha}, \alpha \in \Go \}$
is a Markov process, or a Markov family of probability measures;
see Proposition \ref{strong Markov property}
below. For a measure $\mu$ on $\Go$,
 we define
 
 \[
  P_{\mu} = \int P_{\alpha} \mu (d \alpha).
 \]

We denote by $E_{\mu}$ the expectation under $P_{\mu}$.

\end{rmk}

 \begin{comment}

 \textbf{Check if the following remark needed elsewhere in the text}
 \begin{rmk}\label{3rmk2}
Let $b_1,d_1$ be another pair of birth and death rates, 
satisfying all conditions imposed on 
$b$ and $d$. 
Consider a unique solution $(\tilde{\eta} _t)$ of \eqref{se} with rates 
$b_1, d_1$ instead of
$b,d$, but with the same initial condition $\eta _0$ and all the
other underlying structures. \textit{If for all }
$\zeta \in D$, \textit{where}
$D \in \mathscr{B}( \Gamma _0 (\R^\d))$ , 
$b_1 (\cdot ,\zeta) \equiv b (\cdot ,\zeta) $, $d_1 (\cdot ,\zeta) \equiv d (\cdot ,\zeta) $,
 \textit{then} 
$\tilde{\eta} _t = \eta _t$ \textit{\ for all } $t\leq \inf \{s \geq 0: \eta_s \notin D\} = 
\inf \{s \geq 0: \tilde{\eta}_s \notin D\} $.
This may be proven in the same way as the theorem above.

\end{rmk}
\end{comment}

\begin{rmk}\label{Functional dependence}
 We solved equation \eqref{se} 
 $\omega$-wisely. 
  We can deduce from the proof 
 of Proposition \ref{ex un} that $\theta _n$ and $z_n$ are measurable 
 functions of $\eta _0$ and $N$
 in the sense that, for example, 
 $\theta_1 = F _1 (\eta _0, N)$ a.s.
 for a measurable  function 
 $F_1: \Go
 \times \Gamma( \R _+ \times \R ^\d \times \R _+ \times \gtwo) \to \R _+$.
 As a consequence,
 there is a functional 
 dependence of the solution
 process and the ``input'':
 the process $(\eta _t)_{t\geq 0}$
 is some function of 
 $\eta _0$ and $N$.

\end{rmk}

The following corollary is a consequence of Proposition 
\ref{ex un} and Remark \ref{Functional dependence}.

\begin{cor} \label{flabbergast}
Joint uniqueness in law holds for equation \eqref{se} with initial
distribution $\nu$ satisfying

\[
 \int _{\Go} |\gamma| \nu (d \gamma) < \infty.
\]

\end{cor}

As usually, the Markov property of a solution follows from uniqueness.

 \begin{prop}\label{strong Markov property} \textbf{(The strong Markov property)}
 	Let $b$ and $d$ be time-homogenious. 
 The unique solution $(\eta _t) _{t\in [0,\infty )}$ of \eqref{se} is a strong Markov process
 in the following sense.
 	Let $\tau$ be an a.s. finite $(\mathscr{S} _t, t\geq 0)$-stopping time
 	such that $E |\eta _{\tau}| < \infty$.
 Then 
 \begin{equation}\label{flip out}
 P \{(\eta _{\tau + t}, t \geq 0) \in \mathcal{D} \}  = E P_{\eta _\tau}(\mathcal{D}), 
  \ \ \ \mathcal{D} \in \mathscr{B}(D_{\Go}[0,\infty )).
 \end{equation}
 
 Furthermore, 
 for any $\mathscr{D} \in \mathscr{B}(D _{\Go} [0, \infty))$,
 \begin{equation}\label{flip out2}
 P\{ (\eta _{\tau + t}, t \geq 0) \in \mathscr{D}
 \mid \mathscr{S} _{\tau}\} = 
 P\{ (\eta _{\tau + t}, t \geq 0) \in \mathscr{D}
 \mid \eta _{\tau}\};
 \end{equation}
 that is,
 given $\eta _{\tau}$,  $(\eta _{\tau + t}, t \geq 0)$
 is conditionally independent  of $(\mathscr{S} _t, t\geq 0)$.
 
\end{prop}

\textbf{Proof}. 
For $t \geq 0$
\begin{align} \label{se tau}
\eta _{\tau + t} (B) = \int\limits _{(\tau,\tau + t] \times B \times [0, \infty ) \times \gtwo }
& I _{ [0,b(x,\eta _{s-} )] } (u) I \Bigg\{ \int \limits _{\substack{r \in (s, \tau + t ],
		\\
		v \geq 0}} I_{[0,d(x, \eta _{r-})]}(v)\gamma (dr,dv) =0 \Bigg\}
N(ds,dx,du,d\gamma) \notag \\
+ \int\limits _{B \times  \gtwo } 
I \Bigg\{ \int \limits _{\substack{r \in (\tau, \tau + t ],
		\\
		v \geq 0}} & I_{[0, d(x, \eta _{r-})]}(v)\gamma (dr,dv) =0  \Bigg\} \hat {\eta} _{\tau} (dx, d \gamma ) 
+ \eta _{\tau} (B),  \quad
t \geq 0.
\end{align}
where 
$\hat {\eta} _{\tau}  = \sum\limits _{ \substack{q \in \mathcal{I} \cup \mathcal {J}, \\
												x_q \in \eta _{\tau}	} } (x_q, \gamma _q)$.
 Here we need the strong Markov property of the driving process 
as given in Proposition \ref{Strong MP}.
Note that \eqref{se tau}
can be considered as an equation of the type \eqref{se}
with the unique solution 
is $(\eta _{\tau + t}) _{t\in [0, \infty )}$.
From Proposition \ref{ex un}, Remark \ref{Functional dependence}, and Corollary \ref{flabbergast}
we get \eqref{flip out}. The conditional independence 
\eqref{flip out2} follows from Remark \ref{Functional dependence}.
\qed

\textbf{Proof of Theorem \ref{core thm}.}
 The theorem 
is a 
consequence of Proposition \ref{ex un}, Remark \ref{encompass}, and
Proposition \ref{strong Markov property}. In particular, 
the Markov property of $\{P_{\alpha}, \alpha \in \Go \}$
follows from Corollary \ref{flabbergast}.
\qed

Let $N_1$ be the image of $N$ under the projection
\[
(s,x,u ,\gamma) \mapsto (s,x,u).
\]
The process $N_1$ is a Poisson point process on $\R _+ \times \R ^\d \times \R _+$ with
intensity measure $dsdxdu$.

\textbf{Proof of Proposition} \ref{poignant}.
We have
\[
 \eta _t (B) = \int\limits _{(0,t] \times B \times [0, \infty ) \times \gtwo }
 I _{ [0,b(x,s,\eta _{s-} )] } (u) I \Bigg\{ \int \limits _{\substack{r \in (s,  t ],
 		\\
 		v \geq 0}} I_{[0,d(x,r, \eta _{r-})]}(v)\gamma (dr,dv) =0 \Bigg\}
N(ds,dx,du,d\gamma)
\]
\[
+ \int\limits _{B \times  \gtwo } 
I \Bigg\{ \int \limits _{\substack{r \in (0,  t ],
		\\
		v \geq 0}} I_{[0, d(x,r, \eta _{r-})]}(v)\gamma (dr,dv) =0 \Bigg\} \hat {\eta} _0 (dx, d \gamma ) 
\]
\[
 = \int\limits _{(0,t] \times B \times [0, \infty )  }
 I _{ [0,b(x,s,\eta _{s-} )] } (u) 
N_1(ds,dx,du) + \eta _0 (B)
\]
\[
 - \sum\limits _{q \in \mathcal{I}\cup \mathcal{J}} \int\limits _{(0,t] \times [0, \infty ) } 
 I\{ x_q \in \eta _{r-} \}
 I_{[0, d(x_q,r, \eta _{r-})]}(v) \gamma _q (dr,dv).
\]
Recall that $\eta \cup {x}$
	and $\eta \setminus {x}$ are shorthands for 
	$\eta \cup \{ x \}$ and
	$\eta \setminus \{ x \}$, respectively.
  By Ito's formula (\cite[Chapter 2, Theorem 5.1]{IkedaWat}) for $F \in \mathscr{C} _b$ 
  \begin{align*}
 F(\eta _t) - F(\eta _0)
 = &  \sum\limits _{s \leq t} (F(\eta _s) - F(\eta _{s-})) 
\\
 =&
 \int\limits _{(0,t] \times \mathbf{B}(\mathbf{o}_d, R _F ) \times [0, \infty ) }
 I _{ [0,b(x,s,\eta _{s-} )] } (u)
 \big\{ F(\eta _{s-} \cup x  ) - F(\eta _{s-}) \big\}
N_1(ds,dx,du)
\\
& + \sum\limits _{q \in \mathcal{I}\cup \mathcal{J}} \int\limits _{(0,t] \times [0, \infty ) } 
 I\{ x_q \in \eta _{r-} \}
 I_{[0, d(x_q,r, \eta _{r-})]}(v)
 \big\{ F(\eta _{r-} \setminus x  ) - F(\eta _{r-}) \big\}
 \gamma _q (dr,dv).
  \end{align*}

We can write
\begin{multline*}
 \int\limits _{(0,t] \times \mathbf{B}(\mathbf{o}_d, R _F ) \times [0, \infty ) }
 I _{ [0,b(x,s,\eta _{s-} )] } (u)
 \big\{ F(\eta _{s-} \cup x  ) - F(\eta _{s-}) \big\}
N_1(ds,dx,du) 
\\
 = \int\limits _{(0,t] \times \mathbf{B}(\mathbf{o}_d, R _F )  } b(x,s,\eta _{s-} )
  \big\{ F(\eta _{s-} \cup x  ) - F(\eta _{s-}) \big\} dx ds
\\
 +  \int\limits _{(0,t] \times \mathbf{B}(\mathbf{o}_d, R _F ) \times [0, \infty ) }
 I _{ [0,b(x,s,\eta _{s-} )] } (u) 
  \big\{ F(\eta _{s-} \cup x  ) - F(\eta _{s-}) \big\}
\tilde N _1 (ds,dx,du),
\end{multline*}
where $\tilde N = N - dsdxdu$. Since $F \in \mathscr{C}_b$,
the process
\[
 \int\limits _{(0,t] \times \mathbf{B}(\mathbf{o}_d, R _F ) \times [0, \infty ) }
 I _{ [0,b(x,s,\eta _{s-} )] } (u) 
  \big\{ F(\eta _{s-} \cup x  ) - F(\eta _{s-}) \big\}
\tilde N _1 (ds,dx,du)
\]
is a martingale by item (vi) of Definition \ref{weak solution},
see e.g.  \cite[Section 3 of Chapter 2]{IkedaWat}.
Similarly, 
\[
\sum\limits _{q \in \mathcal{I}\cup \mathcal{J}} \int\limits _{(0,t] \times [0, \infty ) } 
 I\{ x_q \in \eta _{r-} \}
 I_{[0, d(x _q,r, \eta _{r-})]}(v)
 \big\{ F(\eta _{r-} \setminus x  ) - F(\eta _{r-}) \big\}
 \gamma _q (dr,dv)
\]
can be decomposed into a sum of
\[
  \int\limits _{(0,t] } 
 \sum\limits _{x \in \eta _{r-}}
d(x,r, \eta _{r-})]
 \big\{ F(\eta _{r-} \setminus x  ) - F(\eta _{r-}) \big\}
 dr
\]
and a martingale. The desired statement follows. \qed

 \textbf{Proof of Proposition \ref{couple}}.
 Let $\tau _1, \tau _2, ...$ be consecutive jump moments
 of the process $(\xi ^{(1)}_t , \xi ^{(2)}_t)$.
 We will show by induction
that each moment of birth for $(\xi ^{(1)}_t)_{t \in [0,\infty )}$ is
a moment of birth for $(\xi ^{(2)}_t)_{t \in [0,\infty )}$ too, and 
each moment of death for $(\xi ^{(2)}_t)_{t \in [0,\infty )}$
is a moment of death for $(\xi ^{(1)}_t)_{t \in [0,\infty )}$
if the dying particle is in $(\xi ^{(1)}_t)_{t \in [0,\infty )}$.
Moreover, in both cases
the birth or the death occurs at exactly the same place. Here a moment
of birth is a random time at which a new particle appears, a
moment of death is a random time at which an existing particle disappears
from the configuration. The statement formulated
here is in fact equivalent
to \eqref{culprit}.

 Here we deal only with
the base case, the induction step is done in the same way.
We have nothing to show 
if $\tau _1$ is a moment of a birth of $(\xi ^{(2)}_t)_{t \in [0,\infty )}$
or a moment of death of $(\xi ^{(1)}_t)_{t \in [0,\infty )}$.
Assume that a new particle is born for  $(\xi ^{(1)}_t)_{t \in [0,\infty )}$
at $\tau _1$,

$$
\xi ^{(1)}_{\tau _1} \setminus \xi ^{(1)}_{\tau _1 - } =\{x_1 \}.
$$
The process $(\xi ^{(1)})_{t \in [0,\infty )}$ satisfies \eqref{se}, therefore a.s.
$N_1(\{x \} \times \{ \tau _1 \} \times [0,b _1(x_1, {\tau _1}, \xi ^{(1)} _{\tau _1-} )]) = 1$.
Since 
\begin{equation} \label{sulky}
\xi ^{(1)} _{\tau _1-} = \xi ^{(1)} _0 \subset \xi ^{(2)} _0 = \xi ^{(2)} _{\tau _1-},
\end{equation}
by \eqref{culpable} we have $b _1(x_1, {\tau _1}, \xi ^{(1)} _{\tau _1-} )
\subset b _2(x_1, {\tau _1}, \xi ^{(2)} _{\tau _1-} )$, and hence
$$
N(\{x \} \times \{ \tau _1 \} \times [0,b _2(x_1, {\tau _1}, \xi ^{(2)} _{\tau _1-} )] \times \gtwo) = 1,
$$
hence 
$$
\xi ^{(2)}_{\tau _1} \setminus \xi ^{(2)}_{\tau _1 - } =\{x_1 \}.
$$

Now let $\tau _1$ be a moment of 
death for $(\xi ^{(2)}_t)_{t \in [0,\infty )}$,
and let $ \xi ^{(2)}_{\tau _1-} \setminus \xi ^{(2)}_{\tau _1  } =\{x_q \}$
for some $q \in \mathcal{I} \cup \mathcal{J}$ (such a $q$ always exists because of \eqref{faucet},
and is unique).
If $x_q \notin \xi ^{(1)}_{\tau _1-}$, we have nothing to prove. 
Hence we also assume $x_q \notin \xi ^{(1)}_{\tau _1-}$.
 We have a.s. $\gamma _q (\{\tau_1 \} \times [0, d_2(x_q, \tau _1 ,  \xi ^{(2)}_{\tau _1-}  )  ]) = 1$.
By \eqref{sulky} and \eqref{culpable2}, $d_1(x_q, \tau _1 ,  \xi ^{(1)}_{\tau _1-}  ) \geq d_2(x_q, \tau _1 ,  \xi ^{(2)}_{\tau _1-}  )$,
hence
 $$\gamma _q (\{\tau_1 \} \times [0, d_1(x_q, \tau _1 ,  \xi ^{(1)}_{\tau _1-}  )  ]) = 1.$$
 It follows that $ \xi ^{(1)}_{\tau _1-} \setminus \xi ^{(1)}_{\tau _1  } =\{x_q \}$.

\qed

 \section{Aggregation model: proofs}\label{sec4}

The main idea behind our analysis in this section 
is to couple the process $(\eta _t)_{t\geq 0}$
with another 
birth-and-death process,
to which we can apply Lemma \ref{lumpability}.

To do so,
let us introduce another pair of the birth and death rates, 
$b_1, d_1 $, and an initial condition 
$\xi _0 = \eta _0 \cap \Lambda$, such that 
$b_1(x, \eta) = d_1 (x, \eta) = 0 $ for $x \notin \Lambda$, $d_1 (x, \eta) =  a^{-|\eta|} $ for
$x \in \Lambda$, $b_1(x, \eta) \leq b(x, \eta)$ for all $x, \eta$, and
for some constant $c>0$
$$
\int\limits _{\Lambda} b_1(x,\eta) dx = c|\eta \cap \Lambda|, \quad \eta \in \Go.
$$
It follows from \eqref{b aggr} that 
there exists a function $b_1$ satisfying these assumptions.

Functions $b_1$, $d_1$ satisfy conditions of 
Theorem \ref{core thm}.
Furthermore, the conditions
of 
Proposition \ref{couple} are satisfied here: 
for $\eta^1, \eta^2 \in \Go$, $\eta^1 \subset \eta^2$
we have 
\[
 b_1(x,\eta^1)\leq b(x,\eta^1) \leq b(x,\eta^2)
\]
as well as
\[
 d_1(x,\eta^1)\geq d(x, \eta^1) \geq d(x,\eta^2).
\]
 Denote by $(\xi _t)_{t \geq 0}$ the unique solution 
 of \eqref{se} with the birth and death rates $b_1, d_1 $
and initial condition $\xi _0$.
By Proposition \ref{couple}, $\xi _t \subset \eta _t$ hold a.s. for all
$t \geq 0$.

  In this section 
we will work on the canonical probability space
 $$ 
 \big(D_{\Go}[0,\infty) \times D_{\Go}[0,\infty), \mathscr{B}(D_{\Go}[0,\infty) \times D_{\Go}[0,\infty)),
 P_{\alpha} \big),
 $$
where
$P_{\alpha}$ is the push-forward of 
the measure $P$ under 
 \[
 \Omega \ni \omega \mapsto (\eta(\alpha, \cdot),(\xi(\alpha, \cdot)) 
 \in D_{\Go}[0,\infty) \times D_{\Go}[0,\infty).
\]

 Consider the embedded Markov chain of the process $(\xi_t)_{t \geq 0}$,
$Y_k := \xi_{\tau _k}$, where $\tau _k$ are the moments of jumps of $(\xi _t)$. 
It turns out that the process 
$u = \{ u_k \}_{k \in \N} $, where $u_k := |Y_k|$, is a Markov chain too.
Indeed,
the equality
\[
P_{\alpha _1}\{ |Y_1| = k  \} =
P_{\alpha _2}\{ |Y_1| =k  \} , \ \ \ k\in \N , \alpha \in \Go.
\]
holds when $|\alpha _1 \cap \Lambda|=|\alpha _2 \cap \Lambda|$, since both sides are equal to
\[
 \left\{
  \begin{array}{l l}
    \frac{c}{c+a^{-|\alpha _1 \cap \Lambda|}} & \quad \text{if} \quad k=|\alpha _1 \cap \Lambda|+1, \\
    \frac{a^{-|\alpha _1 \cap \Lambda|}}{c+a^{-|\alpha _1 \cap \Lambda|}} & \quad \text{if} \quad k=|\alpha _1 \cap \Lambda|-1 ,\\
     0 & \quad \text{in other cases.}
  \end{array} \right. 
\]
Therefore, 
Lemma \ref{lumpability} is applicable here,
with $f(\cdot) = |\cdot |$.

\textbf{Proof of Proposition} \ref{extinction}. Having in mind the inclusion $\xi _t \subset \eta _t $ ($P_{\alpha}$-a.s.), 
we will prove this proposition for $(\xi _t)$.

It follows from \eqref{transition probabilities}
that the transition probabilities for the Markov chain $\{ u_k \}_{k \in \Z _+} $
are given by
\begin{equation} \label{p_ij}
 p_{i,j} =  P_{\alpha}\{ u_k =j \mid u_{k-1} = i \} = \left\{
  \begin{array}{l l}
    \frac{c}{c+a^{-i}} & \quad \text{if} \quad j=i+1, \\
    \frac{a^{-i}}{c+a^{-i}} & \quad \text{if} \quad j=i-1 ,\\
     0 & \quad \text{in other cases,}
  \end{array} \right. 
\end{equation}
for $i \in \N , j \in \Z _+$,
  and $p_{0,j} = I _ {\{ j =0 \} } $.

Since the zero is a trap and it is accessible
from all other states, 
there are no recurrent states except zero,
and the process $u$ has only two possible
types of behavior on infinity:
\[
P _{\alpha} \{ \exists l \in \N \text{\ s.t.\ } u_l= 0 \text{\ or\ }
\lim\limits _{m \to \infty} u_m=\infty \} = 1.
\]

We will now use 
properties of countable state space Markov chains,
see, e.g.,
 \cite[$\mathsection$ 12,
chapter 1]{Chung}. Chung considers there Markov chain with a reflecting barrier at 0, but we may still
apply those results,
adapting them correspondingly. Denote
$\varrho _m = \prod\limits _{k=1}^m \frac{p_{k,k-1}}{p_{k,k+1}}$. Then the probability
$P_{\alpha}\{ \exists k \in \N \text{\ s.t.\ } u_k=0  \} $ equals to 1 if and only if
$\sum\limits _{j=1}^{\infty} \varrho _j = \nobreak \infty$, 
whichever initial condition ${\alpha}$, $|{\alpha} \cap \Lambda|>0$, we have.
Moreover, if $\sum\limits _{j=1}^{\infty} \varrho _j < \nobreak \infty$
and $P _{\alpha} \{ u_0=q \} =1$ (or, equivalently,
$|\alpha \cap \Lambda| =q$, $q \in \N$), then
$p_q := P _{\alpha} \{ \exists k \in \N \text{\ s.t.\ } u_k=0  \} =
\frac {\sum\limits _{j=q} ^{\infty} \varrho _j } {1+ \sum\limits _{j=1}^{\infty} \varrho _j}$. From \eqref{p_ij}
we see that in our case  $\varrho _j =c^{-j} a^{-\frac {j(j+1)}{2}} $, and

\begin{equation} \label{probability of extinction}
p_q = \frac {\sum\limits _{j=q} ^{\infty} c^{-j} a^{-\frac {j(j+1)}{2}} }
{1+ \sum\limits _{j=1}^{\infty} c^{-j} a^{-\frac {j(j+1)}{2}}} \leq
\frac {\sum\limits _{j=q} ^{\infty} c^{-j} a^{-\frac {j^2}{2}} }
{1+ \sum\limits _{j=1}^{\infty} c^{-j} a^{-\frac {j^2}{2}}} .
\end{equation}

Now, for arbitrary $C >1$ choose $q \in \N$ for which
$c^{-1}a^{-\frac{q}{2}} < C^{-1}$.
For $j>q$ we have $c^{-j} a^{-\frac {j^2}{2}} < c^{-j} a^{-\frac {jq}{2}}=(c^{-1} a^{-\frac {-q}{2}})^j<C^{-j}$,
and
\[
\sum\limits _{j=q} ^{\infty} c^{-j} a^{-\frac {j^2}{2}}<
\sum\limits _{j=q} ^{\infty} C^{-j} = \frac{C^{-q}}{1-C^{-1}},
\] 
so 
that the statement of the proposition  for $(\xi _t)_{t \geq 0}$
follows from \eqref{probability of extinction}. 
\qed

Note that for $(\eta _t)$
the events comprising number of particles going to
infinity
and extinction are not exclusive,
in particular not if
$\int _{\Lambda}b(x, \varnothing)dx > 0 $.
 However,
 it holds that
\begin{equation}\label{egregious}
P\bigg( \{ |\xi _t| = 0 \text{\ for sufficiently large $t$  } \} \cup \{
|\xi _t| \to \infty \ , t \to \infty    \} \bigg) =1
\end{equation}
and
\begin{equation}\label{egregious2}
 P\bigg( \{ |\xi _t| = 0 \text{\ for sufficiently large $t$  } \} \cap \{
|\xi _t| \to \infty \ , t \to \infty    \} \bigg) =0.
\end{equation}

The following equality is also taken from
\cite[$\mathsection$ 12,
chapter 1]{Chung}; for $q>s$, $q,s \in \N$,
and all $\beta$ with $|\beta \cap \Lambda| = q$,

$$P_{\beta} \{ \exists k\in \N: |u_k| = s  \} =
\frac {\sum\limits _{j=q} ^{\infty} \varrho _j (s)} {1+ \sum\limits _{j=s+1}^{\infty} \varrho  _j (s)},
$$
where $\varrho _m (s) = \prod\limits _{k=s+1}^m \frac{p_{k,k-1}}{p_{k,k+1}} = 
c^{-(m-s)}a^{ -\frac12 (m-s)(m+s+1)}$;
in our case
\begin{equation}\label{meadow}
 P_{\beta} \{ \exists k\in \N: |u_k| = s  \} =
\frac{\sum\limits _{j=q} ^{\infty} c^{-(j-s)}a^{ -\frac12 (j-s)(j+s+1)}}
{1+ \sum\limits _{j=s+1}^{\infty} c^{-(j-s)}a^{ -\frac12 (j-s)(j+s+1)}} :=c_{q,s}  <1.
\end{equation}

Note that

\begin{equation}\label{gepraegt}
 c_{q+1,1} \to 0, \ \ \ \ q \to \infty
\end{equation}

\textbf{Proof of Proposition} \ref{verfemen}. Let $(X_k)_{k \in \Z _+}$ be the 
embedded chain of $(\eta _t)_{t\geq 0}$. First we will show
that for all $m \in \N$ and $\alpha \in \Go$,

\begin{equation}\label{feisty}
 P_{\alpha}  \{|X_k \cap \Lambda| =m \text{ infinitely often } \}  =0.
\end{equation}

Let $\beta \in \Go$, $|\beta \cap \Lambda| = m$, $m \in \N$
(the case of $m =0$ is similar, and we do not
write it down). 
Denote $\tilde k = \min \{k\in \N : X_k \cap \Lambda \ne X_0 \cap \Lambda \}$.
Since $\xi _t \subset \eta _t $ holds $P_{\beta}$ - a.s.,

\begin{equation}\label{damp air}
\begin{split}
 P_{\beta} \big\{ |X_k \cap \Lambda| > m, \forall k \geq \tilde k \big\} \geq & \
  P_{\beta} \big\{ |Y_k \cap \Lambda| > m, \forall k \geq 1 \big\} 
 \\
  = P_{\beta} \big\{ u_k > m, & \ \forall k \geq 1  \big\}.
\end{split}
\end{equation}

 By \eqref{meadow}, the probability $P_{\beta} \{ u_k > m,  \forall k \geq 1  \}$
 is positive and does
not depend on $\beta$, $|\beta \cap \Lambda| = m$:
\begin{equation}\label{impertinent}
s_m := 
  P_{\beta} \{ u_k > m,  \forall k \geq 1  \} \geq p_{m,m+1} (1- c_{m+1,m}) >0.
\end{equation}

%Moreover, by \eqref{gepraegt}

%\begin{equation}\label{pochen}
%s_m \to 1, \ \ \ \ m \to \infty
%\end{equation}

Define $k^m_{i}$, $i \in \N$, subsequently by
$k^m_{j+1} = \min \{k > k^m_{j} : |X_k \cap \Lambda| =m \text{ and }
\exists \bar k < k : |X_{\bar k} \cap \Lambda|  \ne m \}$, $k^m _0 =0$.
Note that for all $\beta$
\[
 P _{\beta} \bigg\{ \exists n_0: |X_n \cap \Lambda|=m  \text{ for all } n \geq n_0 \bigg\} =0.
\]

 By the strong Markov property,
\begin{equation}\label{travesty}
\begin{gathered}
P_{\alpha}  \bigg\{|X_k \cap \Lambda| = m \text{ infinitely often } \bigg\} \leq
 P_{\alpha}  \bigg\{k^m_{j} < \infty, \forall j \in \N \bigg\} 
\\
=\prod\limits _{j=1}^{\infty} P_{\alpha}  \big\{k^m_{j+1} < \infty \mid k^m_{j} < \infty \big\} = 0,
\end{gathered}
\end{equation}
by  \eqref{damp air} and \eqref{impertinent}.
% and \eqref{pochen}
 Indeed,
if $P_{\alpha}  \{ k^m_{j} < \infty \} >0$,
then

\[
 P_{\alpha}  \{k^m_{j+1} < \infty \mid k^m_{j} < \infty \} = 
 \frac{E_{\alpha} I_{ \{ k^m_{j} < \infty \} } P_{X_{k^m_{j}}}\{ 
 k ^m _1 < \infty
 \} }
 {E_{\alpha} I_{ \{ k^m_{j} < \infty \} }} 
\]
\[
 \leq\frac{E_{\alpha} I_{ \{ k^m_{j} < \infty \} }\big(1- P_{X_{k^m_{j}}}\{ |X_k \cap \Lambda| > m,
 \forall k \geq \tilde k  \}\big) }
 {E_{\alpha} I_{ \{ k^m_{j} < \infty \} }} 
\]
\[
  \leq \frac{E_{\alpha} I_{ \{ k^m_{j} < \infty \} }
  \big(1- P_{X_{k^m_{j}}}\{ u_k > m,  \forall k \geq 1  \}\big) }
 {E_{\alpha} I_{ \{ k^m_{j} < \infty \} }} = 1 - s_m.
\]
Note that  $1 - s_m <1$ does not depend on $j$, hence \eqref{travesty} follows.
Having proved \eqref{feisty}, we observe that

\begin{equation}\label{cruise}
\begin{split}
 \big\{|\eta _t \cap \Lambda| \to \infty \big\} 
   \cup \big\{ \exists & t': \forall t\geq t', |\eta _t \cap \Lambda| = \varnothing \big\} \\
   = 
   \bigg( \bigcup _{m=1 } ^{\infty} \{|X_k \cap \Lambda| =&m \text{ infinitely often} \} \bigg)^c.
   \end{split}
\end{equation}
Note that if for some element
of probability space $\omega \in \Omega $ 
the process $(\eta _t)_{t \geq 0}$ is stuck in a trap $\gamma$, 
$\gamma \cap \Lambda = \varnothing$,
then $\omega$ belongs to the set on the left-hand side of \eqref{cruise} 
and does not belong to the set $\big\{|X_k \cap \Lambda| =m \text{ infinitely often} \big\}$,
$m \in \N$.

The statement of the proposition follows from \eqref{feisty} and 
\eqref{cruise}. 
\qed

\textbf{Proof of Proposition} \ref{finite death}. 
Define $\widetilde{\eta} _t := \eta _t \cap \Lambda$
and let $\widetilde{X} _k = \widetilde{\eta} _{\varsigma _k}$, where $\varsigma _k$
is the ordered sequence of jumps of $(\widetilde{\eta} _t)_{t\geq 0}$.
Of course, the process $\{ \widetilde{\eta} _t \} _{t \geq 0}$
is not Markov in general, and neither is
$\{ \widetilde{X} _k \} _{k \in \N}$. However, for all $\alpha \in \Gamma _0 (\R^\d) $ the inequality

\[
P _{\alpha} \{ |\widetilde{X}_{1}| - |\widetilde{X}_0| = 1   \} \geq
p_ {|\alpha \cap \Lambda| , |\alpha \cap \Lambda| +1} 
\]
holds,
because for every $\zeta \in \Go$, $\zeta \cap \Lambda = m$,
the integral of the birth rate $b(\cdot,\zeta)$
over $\Lambda$ is larger than $c m$,
and the cumulative death rate in $\Lambda$, $\sum\limits _{x \in \zeta \cap \Lambda}d(x,\zeta)$,
 is less than $ m a^{-m}$. 

The probability of the event that absolutely no death occurs is positive, even when the initial configuration
contains only one point inside $\Lambda $:
\[ 
P _{\alpha}  \bigg\{ |\widetilde{\eta} _t| - |\widetilde{\eta} _{t-}| \geq 0 \text{\ for all \ } t \geq 0   \bigg\} =
P _{\alpha} \bigg\{ |\widetilde{X}_{k+1}| - |\widetilde{X}_k| = 1 \text{\ for all \ } k\in \N  \bigg\}  
\]
\[
= \prod\limits _{k \in \N} P _{\alpha} \left\{ |\widetilde{X}_{k+1}| - |\widetilde{X}_k| = 1 \middle |
|\widetilde{X}_{k}| - |\widetilde{X}_{k-1}| =1 , ..., |\widetilde{X}_1| - |\widetilde{X}_0| =1  \right\} 
\]
\[
\geq \prod\limits _{k \in \N} 
\inf\limits _{\substack{ \zeta \in \Gamma ^0 (\R^\d), \\ |\zeta \cap \Lambda|=|\alpha \cap \Lambda| +k }}
P _{\zeta} \{ |\widetilde{X}_{1}| - |\widetilde{X}_0| = 1   \} 
\]
\[
\geq \prod\limits _{i=|\alpha|} ^{\infty} p_ {i, i+1} = \prod\limits _{i=|\alpha|} ^{\infty} \frac{c}{c+a^{-i}} =
\prod\limits _{i=|\alpha|} ^{\infty} \big(1 - \frac{a^{-i}}{c+a^{-i}} \big) >0,
\]
because 
the series $\sum\limits _{i=|\alpha|} ^{\infty}  \frac{a^{-i}}{c+a^{-i}} $ converges. In particular,
$\prod\limits _{i=m} ^{\infty} p_ {i, i+1} \to 1$ as $ m$ goes to $ \infty $. 
Also,
\begin{equation}\label{zugunsten}
 P _{\alpha _n}  \big\{ |\widetilde{\eta} _t| - |\widetilde{\eta} _{t-}| 
\geq 0 \text{\ for all \ } t \geq 0   \big\} \to 1, \ \ \ |\alpha _n \cap \Lambda| \to \infty.
\end{equation}

It is clear only an a.s. finite  number of deaths inside $\Lambda$ occurs  on 
$\{ \exists t': \forall t\geq t', |\eta _t \cap \Lambda| = \varnothing \}$.
By Proposition \ref{verfemen}, it remains to show 
that only an a.s. finite  number of deaths inside 
$\Lambda$ occurs on $\{|\eta _t \cap \Lambda| \to \infty \} = \{|\widetilde{\eta} _t| \to \infty \} $.
Let us introduce the stopping times $\sigma _n = \inf \{s \in \R: |\widetilde{\eta} _s| \geq n  \}$,
which are finite on $\{|\widetilde{\eta} _t| \to \infty \} $.
Only a finite number of events (births and deaths)
occur until arbitrary finite time $P _{\beta}$-a.s.
for all $\beta \in \Go$, hence for  $n \in \N$

\[ 
P _{\alpha} \Bigl( \{ |\widetilde{\eta} _t| - |\widetilde{\eta} _{t-}| \geq 0 \text{\ for all but finitely many \ } t \geq 0   \} 
\cap \{|\widetilde{\eta} _t| \to \infty \} \Bigr)
\]
\[
\geq P _{\alpha} \Bigl( \{ |\widetilde{\eta} _t| - |\widetilde{\eta} _{t-}| \geq 0 \text{\ for all \ } t \geq \sigma _n  \}
\cap \{|\widetilde{\eta} _t| \to \infty \} \Bigr)
\]
\[
=E _{\alpha} \Big[ I_{\{|\widetilde{\eta} _t| \to \infty \}}
P _{\eta _{\sigma _n} }  \big\{ |\widetilde{\eta} _t| - |\widetilde{\eta} _{t-}| 
\geq 0 \text{\ for all \ } t \geq 0  \big\} \Big].
\]
From
$ |\eta _{\sigma _n}| \geq n $
we have by \eqref{zugunsten}
\[
P _{\eta _{\sigma _n} }  \big\{ |\widetilde{\eta} _t| - |\widetilde{\eta} _{t-}|
\geq 0 \text{\ for all \ } t \geq 0  \big\} \to 1,
\ \ \ n \to \infty.
\]
Therefore, 
\[
P _{\alpha} \Bigl( \{ |\widetilde{\eta} _t| - |\widetilde{\eta} _{t-}| \geq 0 \text{\ for all but finitely many \ } t \geq 0   \} 
\cap \{|\widetilde{\eta} _t| \to \infty \} \Bigr)= P _{\alpha} \{|\widetilde{\eta} _t| \to \infty \}. 
\] 
\qed

Proposition \ref{finite death} is also applicable to $(\xi)_{t \geq 0}$, since
$b_1, d_1 $ satisfy all the conditions imposed on $b,d $.

\textbf{Proof of Theorem} \ref{trajectory-wise asymptotics}.
 First we prove the theorem for $(\xi)_{t \geq 0}$: we prove that 
for $P_{\alpha}$-almost all $\omega \in F_1:= 
 \{ \lim\limits _{t \to \infty} |\xi _t \cap \Lambda|=\infty \}$,

\begin{equation} \label{scum}
  \liminf\limits _{t\to \infty} \frac{|\xi _t \cap \Lambda|}{e^{ct}} > 0.
\end{equation}

 Without loss of generality we assume  $u_0 = |\alpha \cap \Lambda| >0 $.
Let $0=\tau _0 < \tau _1 < \tau _2 < ...$ be the moments of jumps
of $(\xi _t)_{t\geq 0}$, so that $\xi _{\tau _k} = Y _{k}$.
We recall that the random variables $u_n = |Y_n|$
constitute a Markov chain by Lemma \ref{lumpability}.
Note that a.s. on $F_1$,  $u_n > 0$ for all $n \in \N$.
Denote $\psi (n) = cn + n a^{-n}$.
Then
\[
\int\limits _{ \Lambda } b_1(x,Y_k) dx + 
\sum\limits _{x \in Y_k }d_1(x,Y_k) =  c|Y_k|+ |Y_k|a^{-|Y_k|} =\psi (u_k).
\]

By Theorem 12.17 in \cite{KallenbergFound}
there exists an independent of $Y$
 sequence of independent unit exponentials
  $\{\gamma _k \}_{k \in \N}$
 such that
$\gamma _k = \psi (u_k)(\tau _{k} - \tau _{k-1})$ a.s. on $\{\tau _{k} < \infty \} \supset F_1$.
 In particular, $\{\gamma _k \}_{k \in \N}$ 
is independent of $\{ u_k\}_{k \in Z_ +}$.

From Proposition \ref{finite death} we know
that only a finite number of deaths inside $\Lambda$ occur a.s.
In particular, 
there exists a positive finite random variable $\mathbf{m}$ such that
the inequalities
\begin{equation}\label{penibel}
 u_0+ n \geq u_n \geq u_0+n-\mathbf{m}(\omega ) ,  \ \ \ n\in \N
\end{equation}
hold
a.s. on $F_1$.

A.s. on $F_1$

\[
 \tau _n  = \sum\limits_{k=1}^{n-1}(\tau _{k+1} - \tau _k) = 
 \sum\limits_{k=1}^{n-1}\frac{\gamma _k}{\psi (u_k)} \geq 
 \sum\limits_{k=1}^{n-1}\frac{\gamma _k}{u_0+ck}.
\]

Due to Kolmogorov's two-series theorem, the series $\sum\limits_{k=1}^{\infty}\frac{\gamma _k}{u_0+ck}$
is  divergent a.s.   (we recall that $E\gamma _k = D \gamma _k =1$). Hence $\tau _n \to \infty$ a.s.

We will show below that a.s. on $F_1$
\begin{equation}\label{bigot}
 c \tau _n \leq   \ln n  +c \tilde \gamma, \ \ \ n \in \N,
\end{equation}
where
$\tilde \gamma $ is some finite a.s. on $F_1$ random variable. Using \eqref{bigot}, 
we obtain
\[
P_{\alpha}(F_1) \leq 
 P_{\alpha} \big\{ |\xi _t| \geq \frac{e^{ct}}{(\mathbf{m}+1)e^{c\tilde \gamma}}, t\geq 0 \big\} = 
 P_{\alpha} \big\{ |\xi _{\tau _n}| \geq \frac{e^{c\tau _{n+1}}}{(\mathbf{m}+1)e^{c\tilde \gamma}}, n \in \N \big\}
\]
\[
 =P_{\alpha} \big\{ u_n \geq \frac {1}{\mathbf{m}+1} e^{c\tau _{n+1} - c\tilde \gamma}, n \in \N \big\}
 = P_{\alpha} \big\{ \ln (u_n) + \ln (\mathbf{m}+1) \geq 
 c\tau _{n+1} - c \tilde \gamma, n \in \N \big\} \leq P_{\alpha}(F_1).
\]
Therefore, a.s. on $F_1$,
$|\xi _t| \geq \frac{e^{ct}}{(\mathbf{m}+1)e^{c\tilde \gamma}} $ for all $t \geq 0$, and hence
 \eqref{scum} holds.

Inequality \eqref{bigot} follows from the a.s. on $F_1$ convergence
of the series

\begin{equation}\label{scam}
 \sum\limits_{k=1}^{\infty} \bigg(\frac{\gamma _k}{\psi (u_k)} - \frac{1}{ck}\bigg).
\end{equation}

To establish the convergence of \eqref{scam}, we note that

\begin{equation}\label{scam1}
 \sum\limits_{k=1}^{\infty} \bigg(\frac{\gamma _k}{\psi (u_k)} - \frac{\gamma _k}{c u_k}\bigg) 
\end{equation}
converges a.s. on $F_1$ by Kolmogorov's two-series theorem:
\[
 -\sum\limits_{k=1}^{\infty} \bigg(\frac{\gamma _k}{\psi (u_k)} - \frac{\gamma _k}{c u_k}\bigg) = 
 \sum\limits_{k=1}^{\infty} \gamma _k \frac{u_k a ^{-u_k}}{cu_k\psi (u_k)}\leq
\frac{1}{c^2} \sum\limits_{k=1}^{\infty} \gamma _k \frac{ a ^{-u_k}}{u_k}
\]
\[
=\frac{1}{c^2} \sum\limits_{k=1}^{\mathbf{m}} + \frac{1}{c^2} \sum\limits_{k=\mathbf{m}+1}^{\infty} 
\leq \frac{1}{c^2} \sum\limits_{k=1}^{\mathbf{m}} \gamma _k \frac{ a ^{-u_k}}{u_k} + 
\frac{1}{c^2} \sum\limits_{j=1}^{\infty} \gamma _k \frac{ a ^{-j}}{j} < \infty.
\]
The series
\begin{equation}\label{scam1.5}
 \sum\limits_{k=1}^{\infty}\bigg(\frac{\gamma _k}{cu_k} - \frac{1}{cu_k}\bigg)
  = \sum\limits_{k=1}^{\infty}\frac{\gamma _k - 1}{cu_k} 
\end{equation}
 too converges a.s. on $F_1$ by Kolmogorov's theorem, \eqref{penibel}, and since $ \{ \gamma _k \}$
is independent of $\{u_k\}$: using conditioning on  $\{u_k\}$ we get
\[
P _\alpha \left( \left\{ \sum\limits_{k=1}^{\infty}\frac{\gamma _k - 1}{cu_k} \text{ converges } \right\} \cap F_1\right)
 = E _\alpha P _\alpha  \left[ \left\{ \sum\limits_{k=1}^{\infty}\frac{\gamma _k - 1}{cu_k} \text{ converges } \right\}  \cap F_1 \Bigg| \{ u_k\} \right]
 \]
 \[
 =
 E _\alpha P _\alpha \big[ F_1  \big| \{ u_k\} \big] = P _\alpha ( F_1).
\]
Finally,  by \eqref{penibel}
\begin{equation}\label{scam2}
 \sum\limits_{k=1}^{\infty} \bigg(\frac{1}{cu_k} - \frac{1}{ck}\bigg)
\end{equation}
also
converges a.s. on $F_1$.

%, and 
%\begin{equation}\label{scam3}
% \sum\limits_{k=1}^{\infty} (\frac{\gamma _k}{k} - \frac{1}{k})
%\end{equation}
%converges, again by Kolmogorov's two series theorem.
The a.s convergence of the series in \eqref{scam} follows from 
the fact that \eqref{scam1}, \eqref{scam1.5}, and \eqref{scam2} 
converge.

We have thus proved that \eqref{scum} holds a.s. on $F_1$. 
To establish the statement of the theorem, 
note that $\tilde \sigma _n = \inf \{t>0: |\eta _t|\geq n  \}$
is finite on $F$ and a.s.
\[
\big\{ \liminf\limits _{t \to \infty} \frac{|\eta _t \cap \Lambda|}{e^{ct}} =0, |\eta _t | \to \infty \big\}
\subset \big\{ \liminf\limits _{t \to \infty} \frac{|\xi _t |}{e^{ct}} =0  \big\}.
\]
It follows from \eqref{egregious} and \eqref{egregious2} that
\[
 P_{\beta} \big\{\liminf\limits _{t \to \infty} \frac{|\xi _t |}{e^{ct}} =0  \big\} =
 P_{\beta} \big\{(\xi _t)_{t \geq 0} \ \text{ goes extinct} \big\}, 
\ \ \ \ \beta \in \Go.
\]

Therefore, by Proposition \ref{extinction} and the 
strong Markov property 

\[
 P_{\alpha} \big\{\liminf\limits _{t \to \infty} \frac{|\eta _t \cap \Lambda|}{e^{ct}} =0, |\eta _t | \to \infty  \big\}
 = E_{\alpha}P_{\eta_{ \tilde \sigma_n}} \big\{ \liminf\limits _{t \to \infty} \frac{|\eta _t \cap \Lambda|}{e^{ct}} =0, |\eta _t | \to \infty  \big\}
\]
\[
 \leq E_{\alpha}P_{\eta_{\tilde \sigma_n}} \big\{ \liminf\limits _{t \to \infty} \frac{|\xi _t |}{e^{ct}} =0  \big\} \leq \tilde C ^{-n},
\]
where
$\tilde C$ is the constant from Proposition \ref{extinction}.
Since $n$ is arbitrary, 
\[
 P_{\alpha}\big\{\liminf\limits _{t \to \infty} \frac{|\eta _t \cap \Lambda|}{e^{ct}} =0, |\eta _t | \to \infty  \big\} =0.
\]
\qed

\textbf{Proof of Corollary} \ref{growth in average}. 
Let us fix a configuration $\alpha$, $ \alpha \cap \Lambda \ne \varnothing $.
We saw in the proof of Theorem \ref{trajectory-wise asymptotics} that
for  $P_\alpha$-almost all $\omega \in F_1 $ we have 
\[
   |\xi _t| \geq \frac{1}{(\mathbf{m} + 1) e^{c\tilde \gamma}}e^{ct},  \ \ \ t\geq 0,
\]
where $\mathbf{m}$ and $\tilde \gamma$ are a.s. finite on $F_1$ random variables.
Let $G_k$ be the set $\{\omega : \frac{1}{(\mathbf{m} + 1)e^{c\tilde \gamma}} \geq \frac 1k \}$, $k \in \N$.
Then $\bigcup\limits _{k \in \N} G_k \supset F_1$, and, since $P _{\alpha} (F_1) >0$,
$$
P _{\alpha} ( G_k \cap F_1) >0
$$
for some $k \in\N$. Hence
\[
E _{\alpha} |\eta _t \cap \Lambda| \geq E _{\alpha} |\xi _t |I_{G_k \cap F_1} 
\geq \frac 1k e^{ct} P_{\alpha} ( G_k \cap F_1 ). \ \ \ 
\]
\qed

 \section{Appendix}\label{sec5}
  
 \subsection{Markovian functions of a Markov chain}

 Let $(S, \mathscr{B}(S))$ be a Polish (state) space. Consider a
 (time-homogeneous) Markov 
 chain on $(S, \mathscr{B}(S))$ as a family of probability measures
 on $S^ \infty$. Specifically,  on the measurable space
 $
 {(\bar \Omega,\mathscr{F}) = (S^ \infty , \mathscr{B}(S ^\infty ))}
 $
 consider a family of probability measures $\{P_s \}_{s \in S}$ 
 such that for the coordinate mappings 
 \begin{align*}
 X_n: \bar \Omega &\rightarrow S, \\
 X_n (s_1,s_2,&...)  = s_n,
 \end{align*}
 the process $X := \{X_n \}_{n \in \Z _+ }$ is a Markov chain
 satisfying
  for all $s \in S$
 $$
 P_s \{X_0 =s \} =1,
 $$
 $$
 P_s \{ X_{n+m_j}\in A_j, j=1,...,l \mid \mathscr{F} _n \} 
 = P_{X_n} \{ X_{m_j} \in A_j, j=1,...,l  \}.
 $$
 Here $A_j \in \mathscr{B} (S)$,
 $m_j \in \N$, $ l \in \N$, $\mathscr{F} _n = \sigma \{ X_1,...,X_n \}$.
 %The expression on the right side should be understood as the
 %evaluation of the function $P_{q} \{ X_{m_j} \in A_j, j=1,...,k_1  \}$ of $q$
 %at  point $q =X_n$. 
 The space $S$ is separable, hence
 there exists a transition probability 
 kernel
 $Q: S \times \mathscr{B} (S) \rightarrow [0,1]$ such that
 $$
 Q(s,A) = P_s \{ X_1 \in A \}, \ \ \ s\in S, \ A \in \mathscr{B} (S).
 $$
 
 Consider a transformation of the chain $X$, $Y_n = f(X_n)$, where
 $f:S\to \R $ is a Borel-measurable function. 
 Here we formulate sufficient 
 conditions for 
 $Y = \{Y_n \}_{n \in \Z _+ }$  to be 
 a Markov chain. 
 A very similar question was discussed by
 Burke and Rosenblatt \cite{BRosenblatt}
 for discrete space Markov chains.
 The following lemma is proven in \cite[Section 4]{shapenodeath} 
 (sadly the formulation in \cite{shapenodeath} contains
  a typo; the  formulation below is taken
  from \cite{bezborodov2020corrigendum}).

 \begin{lem} \label{lumpability}
 	Assume that for any bounded Borel function $h: S\rightarrow S$
 	
 	\begin{equation}\label{insurgent}
 	E_s h (Y_1) =  E_q h (Y_1) \text{\ whenever \ } f(s)= f(q),
 	\end{equation}
 	Then $Y$ is a Markov chain.
 \end{lem}

  Condition \eqref{insurgent} is the
 equality of distributions of $Y_1$
 under two different measures, $P_s$ and $ P_q$.
   Clearly, this result holds for a Markov chain which is
   not necessarily defined on a canonical state space,
   because the property of a process to be a Markov chain
   depends on its distribution only.

\subsection{Strong Markov property of the driving process}\label{app2}

 Let $N$ be compatible with a right-continuous complete filtration
$\{ \mathscr{F}_t \}_{t \geq 0}$, and $\tau$ be a finite a.s.
$\{ \mathscr{F}_t \}_{t \geq 0}$-stopping time.
For $\gamma \in \gtwo$, $\gamma = \sum\limits _{i} \delta _{(s_i,u_i)}$, let 
$\theta _{\tau} \gamma =  \sum\limits _{i: s_i > \tau } \delta _{(s_i - \tau,u_i)}$.
Also, for $\Xi \in \mathscr{B} (\gtwo)$ we define the shift
\[
 \theta _{\tau} \Xi = \{ \gamma \in \gtwo \mid \theta _{\tau} \gamma \in \Xi \}
\]

 Introduce another
point process $\overline N $ on $\R _+ \times \R ^\d \times \R _+ \times \gtwo $,
\begin{equation}\label{N tau}
 \overline N ([0,s] \times U \times \Xi) = 
N ((\tau,\tau + s] \times U \times \theta _{\tau}  \Xi), 
\ \ \ s>0,  \ U \in \mathscr{B}(\R ^\d \times \R _+), \ \Xi \in \mathscr{B} (\gtwo).
\end{equation}

 \begin{prop}\label{Strong MP}
 The process $\overline N$ is a Poisson point process with
 intensity measure $ds \times dx \times du \times \pi$, independent of 
 $\mathscr{F}_{\tau}$.
\end{prop}
 
 \textbf{Proof}. To prove the proposition, it suffices to show that

 (i) for any $b>a>0$, open bounded $U \subset \R ^\d \times \R _+$ 
 and open $\Xi \subset \gtwo$, 
$\overline N ((a,b) \times U \times \Xi)$
 is a Poisson random variable
with mean $(b-a) \times l(U) \times \pi (\Xi)$, 
where $l$ is the Lebesgue measure 
on $\R ^\d \times \R _+$,
and 

 (ii) for any $b_k>a_k>0$, $k=1,...,m$, open bounded
$U_k \subset \R ^\d$ and open 
$ \Xi _k \subset \gtwo$ such that 
$((a_i,b_i) \times U_i \times \Xi _i ) \cap( (a_j,b_j) \times U_j \times \Xi _j ) = \varnothing $,
$i \ne j$, 
 the collection $\{\overline N ((a_k,b_k) \times U_k \times \Xi _k)\}_{k=1,m}$
 is a finite sequence of 
 independent random variables, independent of $\mathscr{F}_{\tau}$.
 
 Let $\tau _n$ be the sequence of $\{ \mathscr{F}_t \}_{t \geq 0}$-stopping times, 
$\tau _n = \frac{k}{2^n}$ on $\{ \tau \in (\frac{k-1}{2^n},\frac{k}{2^n}] \}$, 
$k \in \N$. Then
 $\tau _n \downarrow \tau$ and
$\tau _n - \tau \leq \frac{1}{2^n}$. The
stopping times $\tau _n$ take only countably many values.
Therefore 
the process $N$ satisfies the strong Markov property for $\tau _n$: 
 the processes $\overline N _n$, defined by
 $$
\overline N _n ([0,s] \times U \times \Xi) :=  N ((\tau _n,\tau _n + s] \times U \times \theta _{\tau _n}  \Xi), 
$$
 are Poisson point processes, independent of $\mathscr{F}_{\tau _n}$, with intensity $ds \times dx \times du \times \pi$.
 
 To prove (i), 
note that $\overline N _n ((a,b) \times U \times \Xi ) \to \overline N ((a,b) \times U \times \Xi) $ a.s. 
and all random variables $\overline N _n ((a,b) \times U \times \Xi)$ have the same distribution, 
therefore $\overline N  ((a,b) \times U \times \Xi)$ is a Poisson random variable with mean 
$(b-a) l (U) \pi (\Xi) $.
The random variables $\overline N _n ((a,b) \times U \times \Xi)$
are independent of $\mathscr{F}_{\tau}$, 
hence $\overline N  ((a,b) \times U \times \Xi)$ is independent 
of $\mathscr{F}_{\tau}$, too. Similarly, the other part of  (ii) follows. \qed

	Let us now show that the filtration $(\mathscr{S}_t)$ defined below \eqref{nonchalant} is
	right-continuous.
	 Indeed, as in the proof 
 	of Proposition \ref{Strong MP}, we can  check
 	that $\overline N _a$ is independent of $ {\mathscr{S}}_{a+}$.
 	Since 
 	$ {\mathscr{S}}_{\infty}  = \sigma(\overline N _a) \vee 
 	 {\mathscr{S}}_{a}$,
 	$\sigma(\tilde N _a)$ and ${\mathscr{S}}_{a}$ are 
 	independent
 	and $ {\mathscr{S}}_{a+}
 	\subset   {\mathscr{S}}_{\infty} $,
 	we see that
 	${\mathscr{S}}_{a+} \subset {\mathscr{S}}_{a} $. Thus, 
 	${\mathscr{S}}_{a+} = {\mathscr{S}}_{a} $.

 \section*{Acknowledgement}

 The authors are thankful to Yuri Kondratiev for numerous discussions
 on the subject.  
 VB was supported by the Department of Computer Science at the University of Verona.
 VB acknowledges a partial support of the DFG through the SFB 701 
 (Bielefeld University) and the IRTG (IGK) 1132 ``Stochastics and Real World Models''.

\bibliographystyle{alphaSinus}
\bibliography{Sinus}

\end{document}